\documentclass[DIV=15, bibliography=totoc]{scrartcl}  

\usepackage[a4paper,margin=1.5cm,bottom=2.5cm]{geometry}
\usepackage{amsmath,amssymb}
\usepackage{algorithmicx}
\usepackage{algpseudocode} 
\usepackage{algorithm}
\usepackage{graphicx, amsthm, bm, lipsum, mathtools, mathrsfs}

\usepackage{multicol}
\usepackage{nomencl}
\makenomenclature

\renewcommand{\nompreamble}{\begin{multicols}{2}}
\renewcommand{\nompostamble}{\end{multicols}}
\usepackage{float}
\usepackage{xcolor} 
\definecolor{lightblue}{rgb}{0,0.5,1.0}
\definecolor{linkblue}{rgb}{0,0.1,0.6}
\definecolor{citegreen}{rgb}{0,0.4,0.0}%{0.1,0.5,0.4}%{0.125,0.6,0.5}
\definecolor{linkred}{rgb}{0.8,0,0.005}%{0.6,0,0.1}
\definecolor{mailviolet}{rgb}{0.3,0,0.35}%{0.6,0,0.1}
\definecolor{tumblue}{rgb}{0,0.396,0.741}
\definecolor{darkgreen}{rgb}{0,0.4,0} 
\definecolor{darkbrown}{rgb}{0.5, 0.396, 0.09}

\usepackage[hyperindex,colorlinks=true,linkcolor=linkblue,citecolor=citegreen,urlcolor=mailviolet,filecolor=linkred]{hyperref}

\usepackage{caption, subcaption}
\usepackage{fancyhdr}

\usepackage{siunitx}
\usepackage[square,comma,nonamebreak,compress,numbers]{natbib}

\usepackage{soulpos}
\usepackage{ulem}
\usepackage{authblk} % used for affiliations 
\usepackage{float}

\usepackage[inkscapearea=page]{svg}
\usepackage{amsmath}
\usepackage{amssymb}
\usepackage[noabbrev]{cleveref}
\usepackage{upgreek}
\usepackage{adjustbox}
\usepackage{color}
\usepackage{enumitem}
\usepackage{todonotes}

\usepackage{pgfplots}
\usepackage{pgfplotstable}
\usepackage{filecontents}
\usepackage{tikz}
\usepackage{tikzscale}
\usepackage{tikz-3dplot}

\usetikzlibrary{spy}
\usetikzlibrary{intersections}
\usetikzlibrary{arrows,shapes}
\usetikzlibrary{spy}
\usetikzlibrary{backgrounds}  
\usetikzlibrary{decorations}
\usetikzlibrary{decorations.markings}
\usetikzlibrary{positioning}
\usetikzlibrary{patterns}
\usetikzlibrary{calc} 
\usetikzlibrary{quotes} 
\usetikzlibrary{external} % Is activated with the command '\tikzexternalize'.
\usetikzlibrary{matrix}

\pgfplotscreateplotcyclelist{customCycleList}
{%
  {blue, mark=o},
  {red, mark=square},
  {darkgreen, mark=diamond},
  {cyan, mark=diamond},
  {darkbrown, mark=triangle*},
  {black, mark=pentagon},
}

\pgfplotsset{every axis/.append style= {
    cycle list name=customCycleList,
}}

\usepackage{abstract}
\usepackage{placeins}

\title{Full-waveform inversion via the scaled boundary finite element method}

\author[1]{Alireza Daneshyar\thanks{\href{mailto:alireza.daneshyar@tum.de}{\texttt{alireza.daneshyar@tum.de}}, Corresponding author}$^{,}$}
\author[2]{Stefan Kollmannsberger}

\affil[1]{Chair of Computational Modeling and Simulation, Technical University of Munich, Germany}
\affil[2]{Data Science in Civil Engineering, Bauhaus-Universität Weimar, Germany}

\newcommand{\publicationDate}{\today}

\date{}
\fancypagestyle{plain}{%
\fancyhf{}%
\vspace{-0.0cm}
\fancyfoot[L]{
\vspace{-0.0cm} 
\footnotesize{\textit{Preprint}\hfill
\publicationDate
\\
}
}
}

\crefname{paragraph}{paragraph}{paragraphs}
\Crefname{paragraph}{Paragraph}{Paragraphs}

\begin{document}
\vspace{-1.5cm} 
\normalem \maketitle  
\normalfont\fontsize{11}{13}\selectfont
\vspace{-1.5cm} \hrule 
\section*{Abstract}
    We begin by addressing the time-domain full-waveform inversion using the adjoint method. Next, we derive the scaled boundary semi-weak form of the scalar wave equation in heterogeneous media through the Galerkin method. Unlike conventional scaled boundary finite element formulations, the resulting system incorporates variable density and two additional terms involving its spatial derivative. As a result, the coefficient matrices are no longer constant over the elements but instead depend on the so-called radial coordinate, rendering the common solution methods such as low-frequency expansion and continued fractions---both of which assume that these matrices can be defined solely by boundary values---inapplicable. In this context, we introduce a radial discretization scheme for solving the transient scalar wave equation with spatially varying density within the framework of the scaled boundary finite element method. This approach begins with temporal discretization in the time domain, followed by spatial discretization along the radial coordinate. We employ the Newmark method for the former and a finite differencing scheme for the latter. However, the choices underlying our ansatz are made for demonstration purposes and remain flexible. Following these discretizations, we introduce an algorithmic condensation procedure to compute the dynamic stiffness matrices of the elements on the fly. Therefore, we maneuver around the need to introduce auxiliary unknowns into the global algebraic system. As a result, the optimization problem is structured in a two-level hierarchy. The higher level computes the nodal wavefields on a relatively coarse mesh, while the lower level provides the gradient of the cost function for a finely distributed set of model parameters. We obtain the Fr\'echet kernel by computing the zero-lag cross-correlations of the forward and adjoint wavefields, and solve the minimization problems iteratively by moving downhill on the cost function hypersurface through the limited-memory BFGS algorithm. To establish a benchmark for assessing the proposed formulation, we omit the density variations along the radial coordinate to allow the resulting system to be approached by conventional solution methods such as low-frequency expansion or continued fractions. We deduce that the low-frequency expansion would suffice for this reduced system, as the bottleneck lies in capturing the density variations, not in resolving the propagating waves. Synthetic data for two flawed samples are generated using standard finite element analysis, ensuring that no inverse crime is commited. Using the low-frequency expansion as the benchmark, the inverse problems are solved and the reconstructed images are compared with those obtained from the proposed radial discretization. The numerical results demonstrate the effectiveness and robustness of the new formulation and show that using the simplified differential equation along with the conventional scaled boundary finite element formulation is highly inferior to applying the complete form of the differential equation. This approach effectively decomposes the computational load into independent local problems and a single coupled global system, making the solution method highly parallelizable and enhancing the scalability of the solver in high-performance computing environments. We demonstrate that, with a simple OpenMP implementation using 12 threads on a personal laptop, the new formulation outperforms the existing approach in terms of both the quality of reconstructed images and computation time.

    \vspace{0.25cm}
    \noindent\textit{Keywords:} 
    scaled boundary finite element method; full-waveform inversion; time reversal; adjoint method; back-propagation; Fr\'echet kernel; gradient-based optimization; inverse problems
    \vspace{-0.4cm}

    \section{Introduction}
    Inhomogeneities within a medium cause wave scattering, which can be translated into visual representations that reveal the underlying mechanical properties of the medium. This process involves a nonlinear optimization technique known as Full-Waveform Inversion (FWI), where the misfit between observed and synthetic data is minimized. Each FWI epoch consists of generating mechanical waves from sources, recording the data at receivers, back-propagating the residuals between the observed and synthetic data in reverse time, computing the gradient of a cost function via zero-lag cross-correlation of the forward and back-propagated wavefields, and updating the model parameters. Thus, each epoch involves solving the wave equation several times. Given the multitude of analyses involved, it is crucial to minimize the footprint of the numerical method to maintain computational efficiency. Nevertheless, high-resolution imaging inherently demands expensive computations as the quality of the resulting images depends on the spatial distribution of the model parameters and, indeed, the shortest wavelength that can be reproduced by the numerical model. Hence, a key ingredient to prevent prohibitive computations is the choice of the numerical method. In this regard, we aim at formulating the scaled boundary finite element method---a powerful numerical tool that excels in reproducing high-frequency waves at low costs---for inverse problems in the wave equation. The so-called \textit{semi-discretized} formulation of the method enables using relatively coarse meshes without losing accuracy. However, the current form of the method cannot handle spatially varying mechanical properties across elements, a capability that is crucial for high-resolution imaging. Addressing this limitation involves deriving the partial differential equations governing the scaled boundary finite elements such that those variations are considered, as well as developing a solution method capable of resolving those variations for the desired number of points within each element.
    
    Full-waveform inversion has found applications in many areas, from geophysical tomography~\cite{zhu2023fourier} to acoustic imaging of the human brain~\cite{guasch2020full}. Although the physical phenomenon is represented by a linear partial differential equation---the governing equation for mechanical waves---the inverse problem is essentially nonlinear since the state of the system depends on the unknown parameters. The task is to minimize an objective function, typically chosen as the L2 norm of the waveform misfit~\cite{seidl2016iterative}. Similar to other nonlinear minimization problems, this process iterates by updating the unknowns in the direction that maximizes the reduction of the data error. The direction is determined by evaluating the change in the residual caused by a perturbation of the current assumed model parameters. It involves computing the Fr\'echet derivatives by solving the wave equation once for each perturbation~\cite{pratt1998gauss}. This prohibitive computation is typically avoided using the adjoint method, which allows us to calculate the sensitivity of the objective function with respect to the perturbation of model parameters with an extra simulation~\cite{gauthier1986two, tromp2005seismic, fichtner2006adjoint}. It states that the gradient terms are equivalent to the weighted zero-lag cross-correlations of the forward and adjoint wavefields~\cite{rao2022point}. The former is the incident field originating from the excitation source, while the latter is the back-propagated field resulting from the sum of residuals emanating in reverse time from the receivers~\cite{plessix2006review}. Hence, as previously mentioned, only one additional simulation that propagates the residuals into the domain is required~\cite{seidl2017full}. Having both wavefields available, the gradient terms can be computed at a computational cost negligible compared to that of solving the wave equation. Although the number of model parameters can reach several million for high-resolution imaging, the efficiency can be preserved if the footprint of the wave equation remains minimal~\cite{virieux2009overview}.
    
    This imaging tachnique was initially limited to two-dimensional applications due to high computational efforts (see for example~\cite{gauthier1986two, igel1996waveform}). Advancements in numerical modeling techniques and computational resources have made it possible to perform simulations in complex and realistic three-dimensional heterogeneous media, as demonstrated by Komatitsch and Vilotte~\cite{komatitsch1998spectral}, Komatitsch and Tromp~\cite{komatitsch1999introduction}, and Peter et al.~\cite{peter2011forward}, among others. Advanced FWI techniques are now used in a wide range of scientific fields, with applications extending to areas such as seismic microzonation~\cite{kubina2018adjoint}, tsunamis~\cite{zhou2019adjoint}, helioseismology~\cite{mandal2017finite}, non-destructive testing~\cite{van2021electrochemical}, and medical imaging~\cite{lucka2021high}. The problem is ill-posed in general in the sense that that it has many local minima \cite{rao2021quantitative}. There are several solution strategies to linearize this ill-posed problem, including steepest descent, conjugate gradient, quasi-Newton, Gauss-Newton, and full-Newton methods~\cite{virieux2009overview}. In each case, the approximation is iteratively refined by moving downhill on the cost function hypersurface. However, this descent does not necessarily lead to the global minimum of the cost function. Successive inversions from longer to shorter period data assist in overcoming this problem as smooth variations in the medium correspond to smooth cost functions with fewer local minima~\cite{fichtner2010full}. Frequency-domain approaches offer a straightforward framework for this strategy by limiting the sequential inversions to a few discrete frequencies or overlapping frequency groups in ascending order~\cite{jeong2012full}. They also provide an efficient framework by transforming the transient wave problem into a linear algebraic system. Conversely, time-domain approaches offer windowing strategies that allow for focusing on specific arrivals, reducing the ill-posedness of the inverse problem by isolating different events~\cite{brossier2009seismic}. Beyond the approaches already discussed, numerous studies over the past decades have focused on advancing full-waveform inversion methodologies, including but not limited to, regularization techniques~\cite{menke2015relationship}, acceleration strategies~\cite{thrastarson2020accelerating}, cost reduction~\cite{van2020accelerated}, automatic differentiation~\cite{zhu2021general}, source encoding~\cite{tromp2019source}, neural networks~\cite{herrmann2023use}, multi-resolution approaches~\cite{burchner2023isogeometric}, advanced discretization schemes~\cite{burchner2023immersed}, and statistical methods~\cite{fichtner2021autotuning}. 
    
    Having mentioned only a fraction of the advancements made, it need be emphasized that a robust numerical tool providing accurate synthetics, paired with a good starting model, is indeed the key ingredient in full-waveform inversion~\cite{kohn2012influence, prieux2013building, asnaashari2013regularized, asnaashari2015time, zhu2016building, datta2016estimating, teodor2021challenges}. The scaled boundary finite element method is an excellent choice as it is one of the most powerful tools for dealing with problems emerging in applied science and engineering. It was first introduced by Song and Wolf~\cite{song1997scaled} and has since evolved to tackle a wide range of applications, including elastostatics~\cite{daneshyar2021general, tian2024automatic, coelho2024enhanced}, dynamics problems~\cite{gravenkamp2018scaled, gravenkamp2020mass, gravenkamp2020three, zhang2021massively, daneshyar2021shooting, zhang2022asynchronous, kuhn2024explicit}, fracture and contact mechanics~\cite{bulling2019high, pramod2019adaptive, zhang2020adaptive, ooi2020polygon, ya2021open, birk2024use}, material inelasticity~\cite{eisentrager2020sbfem}, dimension-reduced structures~\cite{wallner2020scaled, li2020efficient}, thermal analysis~\cite{iqbal2021development}, higher-order continua~\cite{daneshyar2023scaled}, and acoustics~\cite{gravenkamp2017efficient, liu2019automatic, lozano2023domain}, to name a few. Similar to the standard finite element method, it divides the computational domain into finite element spaces, except that the subdivision uses arbitrary star-shaped polygons or polyhedra, depending on the number of spatial coordinates. The region defined by a scaled boundary finite element consists of all points identified by scaling its boundaries relative to a scaling center. Doing so, only the boundaries are discretized, while the radial rays originating from the scaling center are treated analytically. By deriving the weak form of the governing equation, the field variables remain analytical along those rays and the problem is represented in a semi-discretized format. The resulting formulation delivers several appealing features, such as satisfying radiation conditions, reproducing stress singularities, and capturing complicated solution fields~\cite{daneshyar2024radial}. Nevertheless, although the high-frequency content of propagating waves can be resolved by relatively large elements through continued fractions~\cite{bazyar2008continued, birk2012improved, chen2014high}, the element properties are determined at the external boundaries, rendering it impossible to account for spatially varying mechanical properties. As a result, the main asset of the method for full-waveform inversion becomes ineffective. To overcome this limitation, it is necessary to derive the partial differential equations governing the scaled boundary finite elements in a way that accounts for these variations, while also developing a solution method capable of resolving these variations at the desired number of points within each element. Here, we aim at addressing the aforementioned issues and subsequently perform full-waveform inversion using the resulting formulation.
    
    The outline of the paper is as follows: Section~\ref{sec:inverse} introduces the inverse problem and details the adjoint method for obtaining the Fr\'echet kernel. Section~\ref{sec:scaled} delves into the scaled boundary discretization, covering the weak formulation, the low-frequency expansion of the reduced equations, and the radial discretization of the derived system. Section~\ref{sec:numerical} discusses the numerical implementation and presents results from two case studies: a specimen with a hole and another with a weak inclusion. Finally, in Section~\ref{sec:conclusion}, conclusions are drawn.
    
    \section{Inverse problem} \label{sec:inverse}
    Full-waveform inversion is a minimization technique used to generate detailed subsurface images from complete waveforms. Its goal is to minimize the waveform misfit
    \begin{equation} \label{eq:misfit}
        E(\bm{x}) = \frac{1}{2}\sum_r{\int_0^T{(d-u)^2\delta(\bm{x}-\bm{x}^r)dt}},
    \end{equation}
    where $d(\bm{x},t)$ is the measured data, $u(\bm{x},t)$ is the synthetic data, $\delta(\bm{x}-\bm{x}^r)$ is the Dirac delta function representing a receiver located at $\bm{x}^r$, and $[0,T]$ is the time interval of interest. The synthetics $u(\bm{x},t)$ solves the scalar wave equation
    \begin{equation} \label{eq:main}
         \mathcal{L}u = \nabla\cdot(\rho c^2 \nabla u) - \rho\partial^2_t{u} = f,
    \end{equation}
    governing the heterogeneous medium $\Omega$ with spatially varying elastic properties. The operator $\mathcal{L}$ in the above denotes the scalar wave equation, $\nabla$ denotes the spatial gradient, $\partial_t$ denotes the temporal derivative, $\rho(\bm{x})$ is the density of the material, $c(\bm{x})$ is the wave speed in the medium, and $f(t)$ is the source term. The homogeneous initial conditions specifying the state of $\Omega$ at time $t=0$ read
    \begin{equation}
        u(\bm{x},0) = 0, \quad
        \partial_t{u}(\bm{x},0) = 0,
    \end{equation}
    and the free surface boundary condition on the boundary $\partial\Omega$ is
    \begin{equation}
        \nabla u(\bm{x},t)\cdot\bm{n} = 0,
    \end{equation}
    where $\bm{n}$ is the outward-pointing unit normal vector on $\partial\Omega$.
    
    A variety of choices exists to address the inverse problem, including scaling the density, scaling the wave speed, or simultaneously adjusting both parameters. Referring to the comprehensive investigation conducted by B{\"u}rchner et al.~\cite{burchner2023immersed} regarding these choices, the density scaling approach outperforms the others, as it allows for the reconstruction of voids in addition to weak inclusions, whereas downscaling the wave speed leads to oscillations and erroneous scattering when dealing with high material contrasts. Note that, given the relation between density $\rho$, wave speed $c$, and bulk modulus $K$ as
    \begin{equation}
        K = \rho c^2,
    \end{equation}
    this assumption coincides with scaling the bulk modulus as well.
    
    \subsection{Adjoint method}
    Let $H$ be a real Hilbert space, where $H=L^2(\Omega)$ is the space of square-integrable functions on a domain $\Omega\subset\mathbb{R}^n$. For scalar functions $u,v:H\to\mathbb{R}$ the inner product over the temporal and spatial domains reads
    \begin{equation} \label{eq:inner_scalar}
        \langle u,v \rangle = \int_t{\int_\Omega{uv \: d\Omega} dt} .
    \end{equation}
    This inner product is a special case of a bilinear form $B:H\times H\to\mathbb{R}$ satisfying
    \begin{equation}
        \begin{gathered}
            B(\alpha u+v,w) = \alpha B(u,w) + B(v,w), \\
            B(u , \alpha v + w) = \alpha B(u,v) + B(u,w).
        \end{gathered}
    \end{equation}
    Analogous to~\eqref{eq:inner_scalar}, the inner product for vector-valued functions $\bm{u},\bm{v}: H\to\mathbb{R}^m$ is defined as
    \begin{equation} \label{eq:inner_vector}
        \langle \bm{u}, \bm{v} \rangle = \int_t{\int_\Omega{\bm{u}\cdot\bm{v} \: d\Omega} dt} = \int_t{\int_\Omega{\bm{u}^\intercal\bm{v} \: d\Omega} dt},
    \end{equation}
    where $\bm{u}\cdot\bm{v}$ is the standard Euclidean dot product of $\bm{u}$ and $\bm{v}$ in $\mathbb{R}^m$, and $\bm{u}^\intercal\bm{v}$ represents the corresponding vector product. Following from this definition, we can simply deduce that
    \begin{equation} \label{eq:transpose}
        \langle \bm{u}, \bm{A}\bm{v} \rangle = \langle \bm{A}^\intercal\bm{u}, \bm{v} \rangle.
    \end{equation}
    Now, let $\mathcal{L}$ be a linear and bounded operator on $H$. Due to its linearity, we have
    \begin{equation}
        \mathcal{L}(\alpha\bm{u} + \beta\bm{v}) = \alpha\mathcal{L}\bm{u} + \beta\mathcal{L}\bm{v},
    \end{equation}
    for all vector-valued functions $\bm{u}$ and $\bm{v}$ in the Hilbert space $H$ and scalars $\alpha$ and $\beta$ in $\mathbb{R}$. By generalizing the inner product in~\eqref{eq:transpose}, we introduce $\mathcal{L}^\dagger$ as the adjoint of the operator $\mathcal{L}$ through the relation
    \begin{equation}
        \langle \bm{u}, \mathcal{L}\bm{v} \rangle = \langle \mathcal{L}^\dagger\bm{u}, \bm{v} \rangle.
    \end{equation}
    Like $\mathcal{L}$, the adjoint operator $\mathcal{L}^\dagger$ is also linear and bounded~\cite{zeidler2013nonlinear}.
    
    Now, let us proceed with the minimization problem. In order to minimize the waveform misfit in~\eqref{eq:misfit} while satisfying the scalar wave equation, we employ the Lagrange multiplier method, resulting in the following cost function~\cite{abreu2024understanding}
    \begin{equation}
        \chi = \int_0^T{\int_\Omega{\Big[ \frac{1}{2}(d-u)^2\delta(\bm{x}-\bm{x}^r) + \lambda(\mathcal{L}u-f) \Big]d\Omega}dt},
    \end{equation}
    where $\lambda(\bm{x},t)$ is the Lagrange multiplier. Using the bilinear form~\eqref{eq:inner_scalar}, the cost function can be rewritten as
    \begin{equation} \label{eq:cost}
        \chi = \frac{1}{2} \langle e,e \rangle + \langle \lambda,\mathcal{L}u-f \rangle,
    \end{equation}
    where
    \begin{equation}
        e = (d-u)\delta(\bm{x}-\bm{x}^r).
    \end{equation}
    
    Minimizing the discrepancy between the observed and synthetic data, while enforcing the elastic wave equation as a constraint, requires satisfying the Karush--Kuhn--Tucker optimality conditions~\cite{boyd2004convex}. To this end, we express primal feasibility as
    \begin{equation}
        \delta_\lambda \chi = 0, 
    \end{equation}
    and stationarity with respect to $u$ and $\rho$ as
    \begin{equation}
        \delta_u \chi = 0, \quad \delta_\rho \chi = 0.
    \end{equation}
    Primal feasibility gives
    \begin{equation}
        \delta_\lambda \chi = \langle \delta\lambda,\mathcal{L}u-f \rangle = \int_0^T{\int_\Omega{\delta\lambda(\mathcal{L}u-f) d\Omega}dt} = 0,
    \end{equation}
    implying that the integrand must vanish for arbitrary variations $\delta\lambda$. As a result,
    \begin{equation}
        \mathcal{L}u-f = 0
    \end{equation}
    must hold, which is essentially the wave equation governing the synthetics $u$.
    
    Next, we consider the stationarity condition with respect to $u$, which requires that the variation $\delta_u \chi$ vanishes. We introduce an arbitrary but infinitesimally small variation $\delta u$ that vanishes on the boundary $\partial\Omega$, meets the free surface boundary condition, and satisfies homogeneous initial conditions at time $t=0$. Substituting this variation into the objective function yields
    \begin{equation} \label{eq:delta_u}
        \delta_u \chi = - \langle \delta u,e \rangle + \langle \lambda,\mathcal{L}\delta u \rangle.
    \end{equation}
    Elaborating on the second bilinear form, we arrive at
    \begin{equation}
        \langle \lambda,\mathcal{L}\delta u \rangle = \int_0^T{\int_\Omega{\lambda(\mathcal{L}\delta u) d\Omega}dt} = \int_0^T{\int_\Omega{\lambda\nabla\cdot(\rho c^2 \nabla \delta u) d\Omega}dt} - \int_0^T{\int_\Omega{\lambda\rho\partial^2_t{\delta u} d\Omega}dt} .
    \end{equation}
    To remove the differential operators on $\delta u$, we need to perform integration by parts twice. Considering the first integral and applying integration by parts in space, we write
    \begin{equation}
        \begin{gathered}
            \int_0^T{\int_\Omega{\lambda\nabla\cdot(\rho c^2 \nabla \delta u) d\Omega}dt} = \int_0^T{\int_\Omega{\delta u\nabla\cdot(\rho c^2\nabla\lambda) d\Omega}dt} - \int_0^T{\int_{\partial\Omega}{\delta u(\rho c^2\nabla\lambda)\cdot d\bm{s}}dt} + \\
            \int_0^T{\int_{\partial\Omega}{\lambda(\rho c^2 \nabla \delta u)\cdot d\bm{s}}dt} ,
        \end{gathered}
    \end{equation}
    where $d\bm{s}$ is the outward-pointing area element on $\partial\Omega$. The boundary integrals vanish due to the boundary conditions on $\delta u$, leaving
    \begin{equation}
        \int_0^T{\int_\Omega{\lambda\nabla\cdot(\rho c^2 \nabla \delta u) d\Omega}dt} = \int_0^T{\int_\Omega{\delta u\nabla\cdot(\rho c^2\nabla\lambda) d\Omega}dt}.
    \end{equation}
    Similarly, by performing integration by parts twice with respect to time on the second integral, we obtain
    \begin{equation}
        \int_0^T{\int_\Omega{\lambda\rho\partial^2_t{\delta u} d\Omega}dt} = \int_0^T{\int_\Omega{\delta u\rho\partial^2_t\lambda d\Omega}dt} - \left. \int_\Omega{\delta u\rho\partial_t\lambda d\Omega}\right|_0^T + \left. \int_\Omega{\lambda\rho\partial_t{\delta u} d\Omega}\right|_0^T.
    \end{equation}
    Now, the second bilinear form in~\eqref{eq:delta_u} is rewritten as
    \begin{equation}
        \langle \lambda,\mathcal{L}\delta u \rangle = \int_0^T{\int_\Omega{\delta u\big[\nabla\cdot(\rho c^2\nabla\lambda) - \rho\partial^2_t\lambda\big] d\Omega}dt} + \left. \int_\Omega{\delta u\rho\partial_t\lambda d\Omega}\right|_0^T - \left. \int_\Omega{\lambda\rho\partial_t{\delta u} d\Omega}\right|_0^T = \langle \mathcal{L}^\dagger\lambda,\delta u \rangle .
    \end{equation}
    Replacing the above into the variation $\delta_u \chi$ in~\eqref{eq:delta_u} and enforcing stationarity with respect to $u$, we obtain
    \begin{equation}
        \delta_u \chi = - \langle \delta u,e \rangle + \langle \mathcal{L}^\dagger\lambda,\delta u \rangle = \langle \mathcal{L}^\dagger\lambda - e,\delta u \rangle = 0.
    \end{equation}
    This equation must hold for any arbitrary variation $\delta u$. Considering that the boundary terms vanish at $t=0$ due to the homogeneous initial condition on $\delta u$, we arrive at the terminal value problem
    \begin{equation}
        \nabla\cdot(\rho c^2\nabla\lambda) - \rho\partial^2_t\lambda = (d-u)\delta(\bm{x}-\bm{x}^r),
    \end{equation}
    subjected to the terminal conditions
    \begin{equation}
        \lambda(\bm{x},T) = 0, \quad
        \partial_t{\lambda}(\bm{x},T) = 0.
    \end{equation}
    Mathematically, we derive the adjoint of the scalar wave equation from the second bilinear form in~\eqref{eq:delta_u}. However, the derived problem lacks physical meaning, as it does not have an initial condition to define the starting state. Moreover, satisfying the terminal conditions is generally impossible, since these conditions alone are insufficient to uniquely determine the solution without an initial state to evolve from.
    
    Alternatively, one can define the Lagrange multiplier in reverse time as $\lambda(\bm{x},T-t)$. Doing so, the second bilinear form in the cost function yields the convolution of its arguments~\cite{tonti1973variational}. As a result, the terminal value problem is transformed into an initial value problem with homogeneous initial conditions, but reversed in time. This can be interpreted as the back-propagation of the error between the measured and synthetic data recorded at the receiver locations through the domain, from $t=T$ to $t=0$.
    
    The final step is to enforce stationarity with respect to $\rho$. By introducing the variation $\delta\rho$ into the cost function in~\eqref{eq:cost}, we have
    \begin{equation}
        \delta_\rho \chi = \langle \lambda,\delta_\rho\mathcal{L}u \rangle = \int_0^T{\int_\Omega{\lambda\nabla\cdot(\delta\rho c^2 \nabla u) d\Omega}dt} - \int_0^T{\int_\Omega{\lambda\delta\rho\partial^2_t u d\Omega}dt} .
    \end{equation}
    Applying integration by parts in space on the first integral and in time on the second integral leads to
    \begin{equation}
        \begin{gathered}
            \delta_\rho \chi = -\int_0^T{\int_\Omega{\nabla\lambda\cdot(\delta\rho c^2 \nabla u) d\Omega}dt} + \int_0^T{\int_{\partial\Omega}{\lambda(\delta\rho c^2 \nabla u)\cdot d\bm{s}}dt} + \\
            \int_0^T{\int_\Omega{\partial_t\lambda\delta\rho\partial_t u d\Omega}dt} - \left. \int_\Omega{\lambda\delta\rho\partial_t u d\Omega} \right|_0^T.
        \end{gathered}
    \end{equation}
    The second and fourth integrals vanishes due to free surface boundary condition and homogeneous initial conditions, respectively. Thus, we obtain
    \begin{equation}
        \delta_\rho \chi = -\int_0^T{\int_\Omega{\delta\rho(c^2\nabla\lambda\cdot\nabla u - \partial_t\lambda\partial_t u) d\Omega}dt} ,
    \end{equation}
    or, equivalently,
    \begin{equation}
        \delta_\rho \chi = \int_\Omega{\delta\rho \mathcal{K}_\rho d\Omega} ,
    \end{equation}
    where 
    \begin{equation}
        \mathcal{K}_\rho = \int_0^T{(\partial_t\lambda\partial_t u - c^2\nabla\lambda\cdot\nabla u) dt} 
    \end{equation}
    represents the Fr\'echet kernel.
    
    With the variation $\delta_\rho \chi$ at hand, the last Karush--Kuhn--Tucker condition, which is stationarity with respect to $\rho$, is satisfied by setting this variation to zero. However, directly solving $\delta_\rho \chi = 0$ is not possible due to the implicit dependence of the forward and adjoint wavefields on $\rho$. Specifically, since $u$ and $\lambda$ must satisfy their own governing equations that include $\rho$, any change in $\rho$ alters $u$ and $\lambda$, which in turn affects the variation $\chi$ with respect to $\rho$. Consequently, $\delta_\rho \chi$ does not yield a closed-form solution, and setting it directly to zero does not provide a straightforward condition for $\rho$.
    
    To overcome this, an iterative gradient-based optimization approach can be adopted. In each iteration, $\rho$ is updated based on the computed Fr\'echet derivative, which, in fact, represents the sensitivity of $\chi$ with respect to $\rho$. After each update to $\rho$, the wavefields $u$ and $\lambda$ are recomputed according to their governing equations to reflect the new value of $\rho$. This process iterates until convergence, defined by the condition that $\delta_\rho \chi$ becomes sufficiently close to zero, indicating that $\chi$ is minimized with respect to $\rho$.
    
    \section{Scaled boundary discretization} \label{sec:scaled}
    \begin{figure}
        \centering
        \includegraphics[scale=1.0]{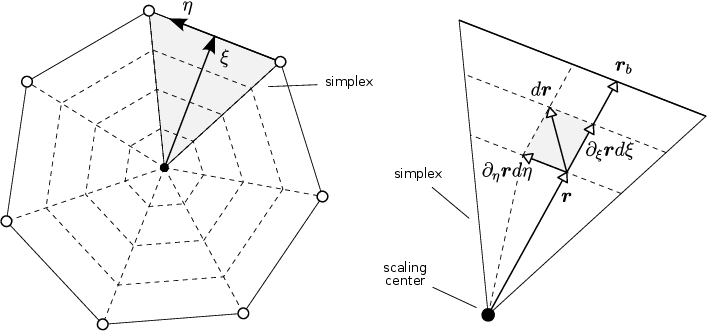}
        \caption{A typical two-dimensional scaled boundary finite element.}
        \label{fig:polygon}
    \end{figure}
    
    We derive the scaled boundary formulation for the scalar wave equation below. We focus on two-dimensional problems here, although formulating three-dimensional cases requires only minor modifications. Following the discussed strategy---assuming variable density with a constant wave speed---the scalar wave equation for two-dimensional setting reads
    \begin{equation} \label{eq:wave}
        \rho c^2 (\partial^2_x{u} + \partial^2_y{u}) + c^2 (\partial_x\rho \partial_x{u} + \partial_y\rho \partial_y{u}) - \rho \partial^2_t{u} = f .
    \end{equation}
    Spatial discretization is carried out within the framework of the scaled boundary finite element method. This method offers flexibility by allowing the use of star-convex polygons and polyhedra with arbitrary numbers of vertices~\cite{song2018scaled}. Figure~\ref{fig:polygon} shows a typical two-dimensional scaled boundary finite element. The polygon is divided into simplices numbering equal to the edges of the element so that their third vertices meet at an interior point known as the scaling center. Subsequently, the transformation between the Cartesian coordinates $(x, y)$ and the local coordinates $(\xi, \eta)$ reads
    \begin{equation} \label{eq:transformation}
        \begin{gathered}
            x(\xi,\eta) = x_0 + \xi (x_b(\eta)-x_0), \\
            y(\xi,\eta) = y_0 + \xi (y_b(\eta)-y_0),
        \end{gathered}
    \end{equation}
    where $(x_0, y_0)$ represent the coordinates of the scaling center, and $(x_b, y_b)$ denote the coordinates of the projection of the point of interest onto the boundary. The ray emanating from the scaling center toward this projection point represents the dimensionless radial coordinate $\xi$. It takes the value of zero at the scaling center and reaches one at the boundary of the element. The local coordinate $\eta$ is the dimensionless circumferential coordinate along which we interpolate through
    \begin{equation}
        \begin{gathered}
            x_b(\eta) = \sum_i{\phi_i(\eta)} x_i  \\
            y_b(\eta) = \sum_i{\phi_i(\eta)} y_i , 
        \end{gathered}
    \end{equation}
    where $\phi_i$ is the shape function of $i$th node and $(x_i, y_i)$ denote its coordinates.
    
    To relate the derivatives in the global coordinate system to those in the local coordinate system, we use the following chain rules
    \begin{equation} \label{eq:chain}
        \begin{gathered}
            \frac{\partial}{\partial\xi}=\frac{\partial}{\partial{x}}\frac{\partial{x}}{\partial\xi}+\frac{\partial}{\partial{y}}\frac{\partial{y}}{\partial\xi}, \\
            \frac{\partial}{\partial\eta}=\frac{\partial}{\partial{x}}\frac{\partial{x}}{\partial\eta}+\frac{\partial}{\partial{y}}\frac{\partial{y}}{\partial\eta}, 
        \end{gathered}
    \end{equation}
    or, equivalently,
    \begin{equation} \label{eq:partial1}
        \begin{bmatrix}
            \partial_\xi  \\
            \partial_\eta 
        \end{bmatrix} =
        \begin{bmatrix}
            \partial_\xi x  & \partial_\xi y \\
            \partial_\eta x & \partial_\eta y \\
        \end{bmatrix}
        \begin{bmatrix}
            \partial_x  \\
            \partial_y 
        \end{bmatrix} ,
    \end{equation}
    by which, the Jacobian matrix is defined as
    \begin{equation}
        \bm J =
        \begin{bmatrix}
            x_b-x_0        & y_b-y_0 \\
            \xi\partial_\eta{x_b} & \xi \partial_\eta{y_b}
        \end{bmatrix} = 
        \begin{bmatrix}
            1 & 0 \\
            0 & \xi
        \end{bmatrix} 
        \begin{bmatrix}
            x_b-x_0    & y_b-y_0 \\
            \partial_\eta{x_b} & \partial_\eta{y_b}
        \end{bmatrix} ,
    \end{equation}
    wherein the coordinate transformation in~\eqref{eq:transformation} is applied. The second square matrix on the right-hand side is the Jacobian matrix at the boundary~\cite{song2018scaled}
    \begin{equation}
        \bm J_b =
        \begin{bmatrix}
            x_b-x_0    & y_b-y_0 \\
            \partial_\eta{x_b} & \partial_\eta{y_b}
        \end{bmatrix} .
    \end{equation}
    Now, the derivatives in the global coordinate system are expressed as
    \begin{equation} \label{eq:partial2}
        \begin{bmatrix}
            \partial_x  \\
            \partial_y 
        \end{bmatrix} = \frac{1}{\xi}
        \begin{bmatrix}
            j_{11} & j_{12} \\
            j_{21} & j_{22}
        \end{bmatrix}
        \begin{bmatrix}
            \xi\partial_\xi  \\
            \partial_\eta
        \end{bmatrix} ,
    \end{equation}
    where
    \begin{equation}
        \begin{aligned}
            & j_{11} = \frac{1}{\det{(\bm J_b)}}\partial_\eta{y_b}, \\
            & j_{12} = \frac{-1}{\det{(\bm J_b)}}(y_b-y_0), \\
            & j_{21} = \frac{-1}{\det{(\bm J_b)}}\partial_\eta{x_b}, \\
            & j_{22} = \frac{1}{\det{(\bm J_b)}}(x_b-x_0).
        \end{aligned}
    \end{equation}
    The derivatives in~\eqref{eq:partial2} can be rewritten as
    \begin{equation} \label{eq:partial3}
        \begin{bmatrix}
            \partial_x  \\
            \partial_y 
        \end{bmatrix} = \bm{b}_1\partial_\xi + \bm{b}_2\frac{1}{\xi}\partial_\eta,
    \end{equation}
    where
    \begin{equation} \label{eq:b1b2}
        \begin{gathered}
            \bm{b}_1 = [j_{11},j_{21}]^\intercal, \\
            \bm{b}_2 = [j_{12},j_{22}]^\intercal.
        \end{gathered}
    \end{equation}
    It can be readily verified from the above that \cite{song1997scaled}
    \begin{equation} \label{eq:substitute}
        \partial_\eta(\det{(\bm J_b)}\bm{b}_2) = -\det{(\bm J_b)}\bm{b}_1.
    \end{equation}
    The above equality will be used to derive the weak form of the problem in the following section.
    
    Referring to the geometrical properties presented in Figure~\ref{fig:polygon}, the position of an arbitrary point inside the polygon is given by means of the vector
    \begin{equation} \label{eq:r}
        \bm{r} = \xi \bm{r}_b,
    \end{equation}
    where $\bm{r}_b$ points to the projection of that point onto the boundary, so that
    \begin{equation} \label{eq:rb}
        \bm{r}_b = 
        \begin{bmatrix}
            x_b - x_0 \\
            y_b - y_0
        \end{bmatrix}.
    \end{equation}
    The total derivative of the position vector $\bm{r}$ reads
    \begin{equation}
        d\bm{r} = \partial_\xi\bm{r}d\xi + \partial_\eta\bm{r}d\eta.
    \end{equation}
    As a result, the infinitesimal elements of length are defined as
    \begin{equation} \label{eq:length_eta}
        |\partial_\eta\bm{r}|d\eta = \xi|\partial_\eta\bm{r}_b|d\eta ,
    \end{equation}
    and
    \begin{equation} \label{eq:length_xi}
        |\partial_\xi\bm{r}|d\xi = |\bm r_b|d\xi .
    \end{equation}
    The infinitesimal area element $d\Omega$ is also computed as
    \begin{equation}
        d\Omega = (\partial_\xi\bm{r}d\xi \times \partial_\eta\bm{r}d\eta)\cdot \bm{e}_z,
    \end{equation}
    where $\bm{e}_z$ is the out-of-plane unit vector. Using~\eqref{eq:r} and~\eqref{eq:rb}, we arrive at
    \begin{equation} 
        d\Omega = \det{(\bm J_b)}\xi d\xi d\eta .
    \end{equation}
    It is worth noting that the unit vectors along $\xi$ and $\eta$ are not necessarily perpendicular, as the infinitesimal area element is actually a parallelogram. Here, we define the unit vectors
    \begin{equation} \label{eq:perpendicular}
        \begin{gathered}
            \bm{n}_{\xi} = \frac{\partial_\eta\bm{r} \times \bm e_z}{|\partial_\eta\bm{r} \times \bm e_z|} = \frac{1}{|\partial_\eta\bm{r}_b|}\det{(\bm J_b)}\bm{b}_1, \\
            \bm{n}_{\eta} = \frac{\bm e_z\times \partial_\xi\bm{r}}{|\bm e_z\times \partial_\xi\bm{r}|} = 
            \frac{1}{|\bm r_b|}\det{(\bm J_b)}\bm{b}_2 , 
        \end{gathered}
    \end{equation}
    which are perpendicular to lines with arbitrary $\xi$ and $\eta$, respectively. We use them later to identify the boundary conditions of the governing equation.
    
    \subsection{weak formulation}
    Following these preliminary derivations, we apply the Galerkin method of weighted residuals to the wave equation in~\eqref{eq:wave}, so that
    \begin{equation} \label{eq:Galerkin}
        \int_\Omega w\rho c^2 (\partial^2_x{u} + \partial^2_y{u}) d\Omega + \int_\Omega w c^2 (\partial_x\rho \partial_x{u} + \partial_y\rho \partial_y{u}) d\Omega - \int_\Omega w\rho\partial^2_t{u} d\Omega = 0,
    \end{equation}
    where $w$ is a weight function. Note that we have omitted the source term on the right-hand side, as it only appears at discrete nodal points and does not contribute to the volume integral. We will incorporate it into the scaled boundary formulation after introducing the internal source vector in the following section.
    
    Using the chain rule in~\eqref{eq:partial3}, the first integral of the above equation is written in the local coordinate system as
    \begin{equation}
        \begin{gathered}
            \int_\Omega w\rho c^2 (\partial^2_x{u} + \partial^2_y{u}) d\Omega  = 
            \int_\xi\int_\eta w\rho c^2 (\bm{b}_1^\intercal \bm{b}_1 \xi \partial^2_\xi{u} + \bm{b}_1^\intercal \bm{b}_2 \partial_\xi\partial_\eta{u} - \bm{b}_1^\intercal \bm{b}_2 \frac{1}{\xi}\partial_\eta{u})\det{(\bm{J}_b)}d\xi d\eta + \\
            \int_\xi\int_\eta w\rho c^2 \bm{b}_2^\intercal\partial_\eta(\bm{b}_1 \partial_\xi{u} + \bm{b}_2\frac{1}{\xi} \partial_\eta{u})\det{(\bm{J}_b)}d\xi d\eta .
        \end{gathered}
    \end{equation}
    By applying integration by parts to the second integral on the right-hand side and utilizing~\eqref{eq:substitute}, we arrive at
    \begin{equation}
        \begin{gathered}
            \int_\Omega w\rho c^2 (\partial^2_x{u} + \partial^2_y{u}) d\Omega =
            \int_\xi\int_\eta w\rho c^2 (\bm{b}_1^\intercal \bm{b}_1 \xi \partial^2_\xi{u} + \bm{b}_1^\intercal \bm{b}_2 \partial_\xi\partial_\eta{u} - \bm{b}_1^\intercal \bm{b}_2 \frac{1}{\xi} \partial_\eta{u})\det{(\bm{J}_b)}d\xi d\eta + \\
            \int_\xi \left. w\rho c^2 \bm{b}_2^\intercal(\bm{b}_1 \partial_\xi{u} + \bm{b}_2\frac{1}{\xi} \partial_\eta{u})\det{(\bm{J}_b)}d\xi \right|_\eta - \\
            \int_\xi\int_\eta (\partial_\eta{w}\rho c^2 \bm{b}_2^\intercal + w\partial_\eta\rho c^2 \bm{b}_2^\intercal - w\rho c^2 \bm{b}_1^\intercal)(\bm{b}_1 \partial_\xi{u} + \bm{b}_2\frac{1}{\xi} \partial_\eta{u})\det{(\bm{J}_b)}d\xi d\eta,
        \end{gathered}
    \end{equation}
    which simplifies to
    \begin{equation}
        \begin{gathered}
            \int_\Omega w\rho c^2 (\partial^2_x{u} + \partial^2_y{u}) d\Omega =
            \int_\xi\int_\eta w\rho c^2 (\bm{b}_1^\intercal \bm{b}_1 \xi \partial^2_\xi{u} + \bm{b}_1^\intercal \bm{b}_1 \partial_\xi{u} + \bm{b}_1^\intercal \bm{b}_2 \partial_\xi\partial_\eta{u})\det{(\bm{J}_b)}d\xi d\eta + \\
            \int_\xi \left. w\rho c^2 \bm{b}_2^\intercal(\bm{b}_1 \partial_\xi{u} + \bm{b}_2\frac{1}{\xi} \partial_\eta{u})\det{(\bm{J}_b)}d\xi \right|_\eta - \\
            \int_\xi\int_\eta (\partial_\eta{w}\rho c^2 \bm{b}_2^\intercal + w\partial_\eta\rho c^2 \bm{b}_2^\intercal)(\bm{b}_1 \partial_\xi{u} + \bm{b}_2\frac{1}{\xi} \partial_\eta{u})\det{(\bm{J}_b)}d\xi d\eta.
        \end{gathered}
    \end{equation}
    Similarly, the second and third integrals in~\eqref{eq:Galerkin} are rewritten in the local coordinate system as
    \begin{equation}
        \begin{gathered}
            \int_\Omega w c^2 (\partial_x\rho \partial_x{u} + \partial_y\rho \partial_y{u}) d\Omega = \\
            \int_\xi\int_\eta w c^2 (\bm{b}_1^\intercal \bm{b}_1\xi\partial_\xi\rho \partial_\xi{u} + \bm{b}_1^\intercal \bm{b}_2 \partial_\xi\rho \partial_\eta{u} + \bm{b}_2^\intercal \bm{b}_1 \partial_\eta\rho \partial_\xi{u} + \bm{b}_2^\intercal \bm{b}_2\frac{1}{\xi}\partial_\eta\rho \partial_\eta{u})\det{(\bm{J}_b)}d\xi d\eta,
        \end{gathered}
    \end{equation}
    and
    \begin{equation}
        \int_\Omega w\rho \partial^2_t{u} d\Omega = \int_\xi\int_\eta w\rho\xi \partial^2_t{u} \det{(\bm{J}_b)}d\xi d\eta,
    \end{equation}
    respectively. Finally, by summing the three integrals in the local coordinate system, the weak formulation in~\eqref{eq:Galerkin} reads
    \begin{equation}
        \begin{gathered}
            \int_\xi\int_\eta w\rho c^2 (\bm{b}_1^\intercal \bm{b}_1 \xi \partial^2_\xi{u} + \bm{b}_1^\intercal \bm{b}_1 \partial_\xi{u} + \bm{b}_1^\intercal \bm{b}_2 \partial_\xi\partial_\eta{u})\det{(\bm{J}_b)}d\xi d\eta - \\
            \int_\xi\int_\eta \partial_\eta{w} \rho c^2 \bm{b}_2^\intercal(\bm{b}_1 \partial_\xi{u} + \bm{b}_2\frac{1}{\xi} \partial_\eta{u})\det{(\bm{J}_b)}d\xi d\eta + 
            \int_\xi\int_\eta w \partial_\xi\rho c^2 \bm{b}_1^\intercal(\bm{b}_1\xi \partial_\xi{u} + \bm{b}_2 \partial_\eta{u})\det{(\bm{J}_b)}d\xi d\eta - \\
            \int_\xi\int_\eta w\rho\xi \partial^2_t{u} \det{(\bm{J}_b)}d\xi d\eta + 
            \int_\xi \left. w\rho c^2 \bm{b}_2^\intercal(\bm{b}_1 \partial_\xi{u} + \bm{b}_2\frac{1}{\xi} \partial_\eta{u})\det{(\bm{J}_b)}d\xi \right|_\eta = 0 .
        \end{gathered}
    \end{equation}
    Due to continuity, the last integral in the above equation cancels out during the assembly of adjacent simplices~\cite{song2018scaled}. Hence, the only remaining part in the weak formulation is
    \begin{equation}
        \begin{gathered}
            \int_\xi\int_\eta w\rho c^2 (\bm{b}_1^\intercal \bm{b}_1 \xi \partial^2_\xi{u} + \bm{b}_1^\intercal \bm{b}_1 \partial_\xi{u} + \bm{b}_1^\intercal \bm{b}_2 \partial_\xi\partial_\eta{u})\det{(\bm{J}_b)}d\xi d\eta - \\
            \int_\xi\int_\eta \partial_\eta{w} \rho c^2 \bm{b}_2^\intercal(\bm{b}_1 \partial_\xi{u} + \bm{b}_2\frac{1}{\xi} \partial_\eta{u})\det{(\bm{J}_b)}d\xi d\eta + 
            \int_\xi\int_\eta w \partial_\xi\rho c^2 \bm{b}_1^\intercal(\bm{b}_1\xi \partial_\xi{u} + \bm{b}_2 \partial_\eta{u})\det{(\bm{J}_b)}d\xi d\eta - \\
            \int_\xi\int_\eta w\rho\xi \partial^2_t{u} \det{(\bm{J}_b)}d\xi d\eta = 0 .
        \end{gathered}
    \end{equation}
    Now, by employing the local coordinate system for the solution field $u$ and the test function $w$ such that
    \begin{equation} \label{eq:interpolation}
        \begin{gathered}
            u(\xi, \eta, t) = \bm{\phi}(\eta) \bm{u}(\xi, t), \\ 
            w(\xi, \eta, t) = \bm{\phi}(\eta) \bm{w}(\xi, t),
        \end{gathered}
    \end{equation}
    the weak form of the problem reads
    \begin{equation} \label{eq:weak}
        \begin{gathered}
            \int_\xi\int_\eta \bm{w}^\intercal\bm{\phi}^\intercal \rho c^2 (\bm{b}_1^\intercal \bm{b}_1 \xi \bm{\phi} \partial^2_\xi\bm{u} + \bm{b}_1^\intercal \bm{b}_1 \bm{\phi} \partial_\xi\bm{u} + \bm{b}_1^\intercal \bm{b}_2 \partial_\eta\bm{\phi} \partial_\xi\bm{u})\det{(\bm{J}_b)}d\xi d\eta - \\
            \int_\xi\int_\eta \bm{w}^\intercal \partial_\eta\bm{\phi}^\intercal\rho c^2 \bm{b}_2^\intercal(\bm{b}_1\bm{\phi} \partial_\xi\bm{u} + \bm{b}_2\frac{1}{\xi} \partial_\eta\bm{\phi}\bm{u})\det{(\bm{J}_b)}d\xi d\eta + \\
            \int_\xi\int_\eta \bm{w}^\intercal\bm{\phi}^\intercal \partial_\xi\rho c^2 \bm{b}_1^\intercal(\bm{b}_1\xi \bm{\phi} \partial_\xi\bm{u} + \bm{b}_2 \partial_\eta\bm{\phi}\bm{u})\det{(\bm{J}_b)}d\xi d\eta - 
            \int_\xi\int_\eta \bm{w}^\intercal\bm{\phi}^\intercal\rho\xi\bm{\phi} \partial^2_t\bm{u} \det{(\bm{J}_b)}d\xi d\eta = 0 .
        \end{gathered}
    \end{equation}
    By defining the distribution of the density field over the simplices through
    \begin{equation}
        \rho(\xi, \eta) = \sum_i{\phi_i(\eta)} \rho_i(\xi) , 
    \end{equation}
    we can introduce
    \begin{equation} \label{eq:Es}
        \begin{gathered}
            \bm{E}_0(\xi) = \int_\eta{\rho c^2 \bm{\phi}^\intercal\bm{b}_1^\intercal \bm{b}_1 \bm{\phi} \det{(\bm J_b)}d\eta}, \\
            \bm{E}_1(\xi) = \int_\eta{\rho c^2 \partial_\eta\bm{\phi}^\intercal\bm{b}_2^\intercal \bm{b}_1 \bm{\phi} \det{(\bm J_b)}d\eta}, \\
            \bm{E}_2(\xi) = \int_\eta{\rho c^2 \partial_\eta\bm{\phi}^\intercal\bm{b}_2^\intercal \bm{b}_2 \partial_\eta\bm{\phi} \det{(\bm J_b)}d\eta}, \\
            \bm{M}_0(\xi) = \int_\eta{\rho \bm{\phi}^\intercal \bm{\phi} \det{(\bm J_b)}d\eta}.
        \end{gathered}
    \end{equation}
    Note that the above square matrices are functions of $\xi$. Substituting~\eqref{eq:Es} into~\eqref{eq:weak} yields
    \begin{equation}
        \int_\xi \bm{w}^\intercal\big( \xi\bm{E}_0 \partial^2_\xi\bm{u} + (\bm{E}_0+\bm{E}_1^\intercal-\bm{E}_1+\xi\partial_\xi\bm{E}_0) \partial_\xi\bm{u} - \frac{1}{\xi}(\bm{E}_2-\xi \partial_\xi\bm{E}_1^\intercal)\bm{u} - \xi\bm{M}_0 \partial^2_t\bm{u} \big)d\xi = 0 ,
    \end{equation}
    which is satisfied for any arbitrary $\bm{w}$ only if
    \begin{equation} \label{eq:ODE}
        \xi^2\bm{E}_0 \partial^2_\xi\bm{u} + \xi(\bm{E}_0+\bm{E}_1^\intercal-\bm{E}_1+\xi\partial_\xi\bm{E}_0) \partial_\xi\bm{u} - (\bm{E}_2-\xi \partial_\xi\bm{E}_1^\intercal)\bm{u} - \xi^2\bm{M}_0 \partial^2_t\bm{u} = \bm{0} .
    \end{equation}
    
    Now, we must extract the internal source vector that corresponds to the above differential equation. To do this, we use the internal flux
    \begin{equation}
        q = \rho c^2\nabla{u}\cdot\bm{n},
    \end{equation}
    which acts on the arbitrary boundary $\Gamma$ having the normal vector $\bm{n}$. Based on the Galerkin method of weighted residuals, for a line with constant $\xi$, we have~\cite{song2018scaled}
    \begin{equation}
        \int_\Gamma{w\rho c^2\nabla{u}\cdot\bm{n}}d\Gamma = \int_\eta{w\rho c^2\bm{b}^\intercal_1(\bm{b}_1\xi\partial_\xi{u} + \bm{b}_2\partial_\eta{u})\det{(\bm J_b)}}d\eta ,
    \end{equation}
    where the chain rule in~\eqref{eq:partial3}, the infinitesimal element of length in~\eqref{eq:length_eta}, and the normal vector $\bm{n}_\xi$ in~\eqref{eq:perpendicular} are used. Introducing the local coordinate system in~\eqref{eq:interpolation} to the above, we arrive at
    \begin{equation}
        \int_\eta \bm{w}^\intercal\bm{\phi}^\intercal \rho c^2 \bm{b}_1^\intercal(\bm{b}_1 \xi \bm{\phi} \partial_\xi\bm{u} + \bm{b}_2 \partial_\eta\bm{\phi} \bm{u})\det{(\bm{J}_b)}d\eta  = \bm{w}^\intercal \bm{q},
    \end{equation}
    where
    \begin{equation} \label{eq:internal}
        \bm{q} = \xi\bm{E}_0 \partial_\xi\bm{u} + \bm{E}_1^\intercal\bm{u}
    \end{equation}
    is known as the internal source vector~\cite{song1997scaled}.
    
    \subsection{Low-frequency expansion}
    One can omit the variation of density along the radial coordinate $\xi$, and only consider its variation along the circumferential coordinate $\eta$. As a result, the differential equation~\eqref{eq:ODE} reduces to
    \begin{equation} \label{eq:simplified}
        \xi^2\bm{E}_0 \partial^2_\xi\bm{u} + \xi(\bm{E}_0+\bm{E}_1^\intercal-\bm{E}_1) \partial_\xi\bm{u} - \bm{E}_2\bm{u} - \xi^2\bm{M}_0 \partial^2_t\bm{u} = \bm{0} .
    \end{equation}
    The resulting differential equation can be approached by conventional solution methods such as low-frequency expansion or continued fractions.
    
    The achievable resolution in full-waveform inversion depends on various factors, such as scattering angles, reflected arrivals, and the distance of the region of interest from sources and receivers. At worst, however, this resolution can be reduced to half the wavelength of the highest available frequency \cite{dessa2003combined, warner2013anisotropic}. Having mentioned that, based on the sensitivity analysis made in the work of Song~\cite{song2009scaled}, three to four continued-fraction terms per wavelength lead to highly accurate results. This means that the low-frequency expansion, equivalent to using one term of the continued-fraction method, can accurately resolve a quarter of a wavelength. As a result, two nodes located symmetrically on opposite sides of the scaling center can roughly provide the lowest possible resolution that full-waveform inversion can capture. Therefore, the bottleneck lies in capturing the density variations, not in resolving the propagating waves. Resultingly, the low-frequency expansion would suffice for the case where the density variations along the radial lines is neglected. Needless to mention that this assumption is inferior to the case where the complete form of the differential equation is applied.
    
    The low-frequency expansion solution of the simplified equation in~\eqref{eq:simplified} is addressed. We first introduce the exponential Fourier transform
    \begin{equation}
            \hat{\mathscr{F}}(\omega) = \int_{-\infty}^{+\infty}{\mathscr{F}(t)\exp(-i\omega t)dt}, \\
    \end{equation}
    where $i$ is the imaginary unit and $\omega$ is the angular frequency. By apply the Fourier transform to the partial differential equation in~\eqref{eq:simplified}, we arrive at the second-order ordinary differential equation system
    \begin{equation}
        \xi^2\bm{E}_0 \partial^2_\xi\bm{u} + \xi(\bm{E}_0+\bm{E}_1^\intercal-\bm{E}_1) \partial_\xi\bm{u} - \bm{E}_2\bm{u} + \omega^2\xi^2\bm{M}_0 \bm{u} = \bm{0} .
    \end{equation}
    Now, by introducing the dynamic stiffness matrix $\bm{S}(\xi,\omega)$, such that
    \begin{equation}
        \bm{q} = \bm{S} \bm{u}, 
    \end{equation}
    this system can be transformed into the first-order nonlinear differential equation
    \begin{equation} \label{eq:nonlinearODE}
        (\bm{S} - \bm{E}_1)\bm{E}_0^{-1}(\bm{S} - \bm{E}_1^\intercal ) - \bm{E}_2 + \omega\partial_\omega\bm{S} + \omega^2\bm{M}_0 = \bm{0} ,
    \end{equation}
    defined on the boundary $\xi=1$. It is known as the scaled boundary finite element equation in dynamic stiffness~\cite{song2008evaluation, song2018scaled}. Following this definition, the low-frequency expansion solution obtained by assuming
    \begin{equation} \label{eq:S}
        \bm{S} \approx \bm{K} - \omega^2 \bm{M},
    \end{equation}
    where $\bm{K}$ and $\bm{M}$ are the stiffness and mass matrices of the element, respectively. Considering the static case, i.e., $\omega=0$, the system reduces to
    \begin{equation}
        (\bm{K} - \bm{E}_1)\bm{E}_0^{-1}(\bm{K} - \bm{E}_1^\intercal ) - \bm{E}_2 = \bm{0} ,
    \end{equation}
    which is an algebraic Riccati equation with standard solution procedures~\cite{song1997scaled}. With the stiffness matrix $\bm{K}$ at hand, substitution of the dynamic stiffness matrix $\bm{S}$ from~\eqref{eq:S} into the differential equation in~\eqref{eq:nonlinearODE} yields
    \begin{equation}
        (\bm{K} - \bm{E}_1)\bm{E}_0^{-1}\bm{M} + \bm{M}\bm{E}_0^{-1}(\bm{K} - \bm{E}^\intercal_1) + 2\bm{M}-\bm{M}_0 = \bm{0},
    \end{equation}
    which is a Lyapunov equation for the mass matrix $\bm{M}$ that can be solved using standard matrix factorization methods. Once the stiffness and mass matrices of the elements are computed, the global algebraic system can be assembled similar to the standard finite element and discretized in the time domain according to the chosen method.
    
    \subsection{Radial discretization}
    The original system of differential equations comprising the density variation requires special treatment as the coefficient matrices vary across the element and cannot be determined solely on the boundary. Additionally, due to the spatially varying density, two extra terms are introduced that comprise the derivatives of coefficient matrices. Thus, we extend our radial discretization method previously presented in~\cite{daneshyar2024radial} for solving the steady-state wave equation in homogeneous elastic media in the frequency domain. We address both the spatially varying density and the transient nature of the equations in this extension.
    
    We proceed with temporal discretization using the Newmark method. However, the derived system of equations can also be discretized in the time domain using other methods. According to the Newmark method, we have
    \begin{equation} \label{eq:Newmark}
        \begin{gathered}
            \partial_t\bm{u}^{t+\Delta t} = \frac{\gamma}{\beta\Delta t}(\bm{u}^{t+\Delta t} - \bm{u}^t) + (1-\frac{\gamma}{\beta})\partial_t\bm{u}^t + \Delta t (1-\frac{\gamma}{2\beta})\partial^2_t\bm{u}^t , \\
            \partial^2_t\bm{u}^{t+\Delta t} = \frac{1}{\beta\Delta t^2}(\bm{u}^{t+\Delta t} - \bm{u}^t) - \frac{1}{\beta\Delta t}\partial_t\bm{u}^t + (1-\frac{1}{2\beta})\partial^2_t\bm{u}^t ,
        \end{gathered}
    \end{equation}
    where the superscript $t$ denotes the current known state, $\Delta t$ is the time increment, $t+\Delta t$ denotes the updated unknown state, and $\beta$ and $\gamma$ are the coefficients pf the method that control the stability and accuracy of the solution. Using the above approximations, the differential equation system in~\eqref{eq:ODE} yields
    \begin{equation} \label{eq:ODE_Newmark}
        \xi^2\bm{E}_0\partial^2_\xi\bm{u}^{t+\Delta t} + \xi(\bm{E}_0+\bm{E}_1^\intercal-\bm{E}_1+\xi\partial_\xi\bm{E}_0)\partial_\xi\bm{u}^{t+\Delta t} - (\bm{E}_2  - \xi\partial_\xi\bm{E}_1^\intercal + \frac{1}{\beta\Delta t^2}\xi_i^2 \bm{M}_0)\bm{u}^{t+\Delta t} = \bm{p}^{t+\Delta t} ,
    \end{equation}
    where
    \begin{equation}
        \bm{p}^{t+\Delta t} = -\xi^2\bm{M}_0 \Big(\frac{1}{\beta\Delta t^2} \bm{u}^t + \frac{1}{\beta\Delta t}\partial_t\bm{u}^t - (1-\frac{1}{2\beta})\partial^2_t\bm{u}^t\Big)
    \end{equation}
    is the residual source vector emerging due to the state of the element at the current known state. Having eliminated the second time derivative, we now proceed with spatial discretization. To this end, we use the second-order central difference approximation so that the spatial derivatives of the vector field $\bm{u}$ at $\xi_i$ are given by
    \begin{equation}
        \begin{gathered}
            \partial_\xi\bm{u}(\xi_i) = \partial_\xi\bm{u}_i + O(h^2), \\
            \partial^2_\xi\bm{u}(\xi_i) = \partial^2_\xi\bm{u}_i + O(h^2), 
        \end{gathered}
    \end{equation}
    where $O(h^2)$ represents the truncation error order of the approximation, $h$ is the spacing between two consecutive points, and
    \begin{equation} \label{eq:FD}
        \begin{gathered}
            \partial_\xi\bm{u}_i = \frac{1}{2h} (-\bm{u}_{i-1} + \bm{u}_{i+1}), \\
            \partial^2_\xi\bm{u}_i = \frac{1}{h^2} (\bm{u}_{i-1} - 2\bm{u}_{i} + \bm{u}_{i+1}).
        \end{gathered}
    \end{equation}
    Using the above approximation, the system in~\eqref{eq:ODE_Newmark} is rewritten as
    \begin{equation} \label{eq:discretized}
        \xi_i^2\bm{E}_0\partial^2_\xi\bm{u}_i^{t+\Delta t} + \xi_i(\bm{E}_0+\bm{E}_1^\intercal-\bm{E}_1+\xi_i\partial_\xi\bm{E}_0)\partial_\xi\bm{u}_i^{t+\Delta t} - (\bm{E}_2 - \xi_i\partial_\xi\bm{E}_1^\intercal + \frac{1}{\beta\Delta t^2}\xi_i^2 \bm{M}_0)\bm{u}_i^{t+\Delta t} = \bm{p}_i^{t+\Delta t} ,
    \end{equation}
    for $i\in\lbrace 0,1,2,\cdots,n \rbrace$. Accordingly, $\xi_0$ and $\xi_n$ denote the scaling center and external boundary of the element, respectively, while the other vectors between them correspond to the interior of the element.
    
    The above relation produces $n+1$ equations which can be written in the following simplified form
    \begin{equation}
        \bm\vartheta_i \bm{u}_{i-1}^{t+\Delta t} + \bm\varphi_i \bm{u}_i^{t+\Delta t} + \bm\psi_i \bm{u}_{i-1}^{t+\Delta t} = \bm p_i^{t+\Delta t},
    \end{equation}
    where
    \begin{equation}
        \begin{gathered}
            \bm\vartheta_i = \frac{1}{h^2}\xi_i^2\bm{E}_0 - \frac{1}{2h}\xi_i (\bm{E}_0 + \bm{E}_1^T - \bm{E}_1 + \xi_i\partial_\xi \bm{E}_0)  , \\
            \bm\varphi_i   = -\frac{2}{h^2}\xi_i^2 \bm{E}_0 - \bm{E}_2 + \xi_i \partial_\xi \bm{E}_1^\intercal - \frac{1}{\beta\Delta t^2}\xi_i^2 \bm M_0 ,\\
            \bm\psi_i      = \frac{1}{h^2}\xi_i^2\bm{E}_0 + \frac{1}{2h}\xi_i (\bm{E}_0 + \bm{E}_1^T - \bm{E}_1 + \xi_i\partial_\xi \bm{E}_0) .
        \end{gathered}
    \end{equation}
    Note that the derivatives $\partial_\xi\bm{E}_0$ and $\partial_\xi\bm{E}_1^\intercal$ are approximated using the central difference approximation similar to those presented in~\eqref{eq:FD}, except for $i=0$ and $i=n$ for which we use forward and backward differencing, respectively. As a result, we have
    \begin{equation}
        \partial_\xi \bm{E}_0(\xi_0) = \frac{1}{h}(-\bm{E}_0(\xi_0)+ \bm{E}_0(\xi_1)), \quad
        \partial_\xi \bm{E}_1(\xi_0) = \frac{1}{h}(-\bm{E}_1(\xi_0)+ \bm{E}_1(\xi_1)),
    \end{equation}
    for $i=0$,
    \begin{equation}
        \partial_\xi \bm{E}_0(\xi_n) = \frac{1}{h}(-\bm{E}_0(\xi_{n-1})+ \bm{E}_0(\xi_n)), \quad
        \partial_\xi \bm{E}_1(\xi_n) = \frac{1}{h}(-\bm{E}_1(\xi_{n-1})+ \bm{E}_1(\xi_n)),
    \end{equation}
    for $i=n$, and
    \begin{equation}
        \partial_\xi \bm{E}_0(\xi_i) = \frac{1}{2h}(-\bm{E}_0(\xi_{i-1})+ \bm{E}_0(\xi_{i+1})), \quad
        \partial_\xi \bm{E}_1(\xi_i) = \frac{1}{2h}(-\bm{E}_1(\xi_{i-1})+ \bm{E}_1(\xi_{i+1})),
    \end{equation}
    otherwise.
    
    In addition to the boundary and interior vectors, each containing the solution field distributed along the radial lines with different $\xi$, two ghost vectors $\bm{u}_{-1}$ and $\bm{u}_{n+1}$, located beyond the boundaries, are introduced. Thus, we use the boundary conditions of the differential equation to define two more equations. By replacing $\xi=0$ and $\xi=1$ in the the internal source vector $\bm{q}$ given in~\eqref{eq:internal}, we arrive at
    \begin{equation}
        \xi_0 \bm{E}_0\partial_\xi\bm{u}_0 + \bm{E}_1^\intercal \bm{u}_0 = \bm{0}
    \end{equation}
    for $\xi=0$, and
    \begin{equation}
        \xi_n \bm{E}_0\partial_\xi\bm{u}_n + \bm{E}_1^\intercal \bm{u}_n = \bm{f}
    \end{equation}
    for $\xi=1$, where $\bm{f}$ is the external source vector acting on the boundary nodes. Note that due to equilibrium, the right-hand side of the boundary condition at $\xi=0$ is a null vector. Now, by applying the central difference approximation~\eqref{eq:FD} on the boundary conditions, two additional equations read
    \begin{equation}
        \begin{gathered}
            \bm\vartheta_{-1} \bm{u}_{-1} + \bm\varphi_{-1} \bm{u}_0 + \bm\psi_{-1}\bm{u}_1 = \bm{0} \quad \text{at} \quad \xi=0, \\
            \bm\vartheta_{n+1} \bm{u}_{n-1} + \bm\varphi_{n+1} \bm{u}_n + \bm\psi_{n+1} \bm{u}_{n+1} = \bm{f} \quad \text{at} \quad \xi=1,
        \end{gathered}
    \end{equation}
    where
    \begin{equation}
        \begin{gathered}
            \bm\vartheta_{-1}  = -\frac{1}{2h}\xi_0 \bm{E}_0, \quad
            \bm\varphi_{-1}    =  \bm{E}_1^\intercal, \quad
            \bm\psi_{-1}       =  \frac{1}{2h}\xi_0 \bm{E}_0, \\
            \bm\vartheta_{n+1} = -\frac{1}{2h}\xi_n \bm{E}_0, \quad
            \bm\varphi_{n+1}   =  \bm{E}_1^\intercal, \quad
            \bm\psi_{n+1}      =  \frac{1}{2h}\xi_n \bm{E}_0.
        \end{gathered}
    \end{equation}
    Finally, we arrive at the following algebraic system of equations
    \begin{equation} \label{eq:system}
        \left[
        \begin{array}{lllllll}
            \bm\vartheta_{-1}&\bm\varphi_{-1} &\bm\psi_{-1}  &                  &                  &                &             \\
            \bm\vartheta_{0} &\bm\varphi_{0}  &\bm\psi_{0}   &                  &                  &                &             \\
                             &\bm\vartheta_{1}&\bm\varphi_{1}&\bm\psi_{1}       &                  &                &             \\
                             &                &\ddots        &\ddots            &\ddots            &                &             \\
                             &                &              &\bm\vartheta_{n-1}&\bm\varphi_{n-1}  &\bm\psi_{n-1}   &             \\
                             &                &              &                  &\bm\vartheta_{n}  &\bm\varphi_{n}  &\bm\psi_{n}  \\
                             &                &              &                  &\bm\vartheta_{n+1}&\bm\varphi_{n+1}&\bm\psi_{n+1}
        \end{array} \right] \left[
        \begin{array}{l}
            \bm{u}_{-1}^{t+\Delta t} \\ \bm{u}_0^{t+\Delta t} \\ \bm{u}_1^{t+\Delta t} \\ \vdots \\ \bm{u}_{n-1}^{t+\Delta t} \\ \bm{u}_{n}^{t+\Delta t} \\ \bm{u}_{n+1}^{t+\Delta t}
        \end{array} \right] = \left[
        \begin{array}{l}
            \bm{0} \\ \bm{p}_0^{t+\Delta t} \\ \bm{p}_1^{t+\Delta t} \\ \vdots \\ \bm{p}_{n-1}^{t+\Delta t} \\ \bm{p}_n^{t+\Delta t} \\ \bm{0}
        \end{array} \right] + \left[
        \begin{array}{l}
            \bm{0} \\ \bm{0} \\ \bm{0} \\ \vdots \\ \bm{0} \\ \bm{0} \\ \bm{f}^{t+\Delta t}
        \end{array}
        \right].
    \end{equation}
    Note that $\bm{p}_{-1}^{t+\Delta t}$ and $\bm{p}_{n+1}^{t+\Delta t}$ are null since the boundary conditions do not involve dynamic terms. Note that $\xi_0$ corresponds to the scaling center at which the system in~\eqref{eq:ODE} is singular. To avoid this singularity, instead of setting $\xi_0$ to zero, we can assign a small value to it. However, this value should not be too small as to cause ill-conditioning. According to the work of Song and Vrcelj~\cite{song2008evaluation}, setting $\xi_0=10^{-3}$ is suitable for most engineering applications.
    
    Using the coefficient matrix in the algebraic system above, we can compute the dynamic stiffness matrix of the element. To this end, we condense the effect of the coefficient matrix on the nodal degrees of freedom located at the $n$th radial point. We can efficiently perform this condensation by treating the coefficient in~\eqref{eq:system} as a tridiagonal matrix. Based on the Gauss elimination method, we first define the pivot matrix
    \begin{equation}
        \bm\varrho_0 = \bm\vartheta_0\bm\vartheta_{-1}^{-1},
    \end{equation}
    by which we eliminate the sub-diagonal term $\bm\vartheta_{0}$ through
    \begin{equation}
        \begin{gathered}
            \bm\vartheta_0 := \bm\vartheta_0 - \bm\varrho_0\bm\vartheta_{-1} = \bm{0}, \\
            \bm\varphi_0   := \bm\varphi_0   - \bm\varrho_0\bm\varphi_{-1}, \\
            \bm\psi_0      := \bm\psi_0      - \bm\varrho_0\bm\psi_{-1}.
        \end{gathered}
    \end{equation}
    Note that the symbol $:=$ is the assignment operator, indicating that the expression does not represent an equation but shows that the value on the left-hand side is being updated.
    
\begin{figure}
        \centering
        \includegraphics[scale=1.0]{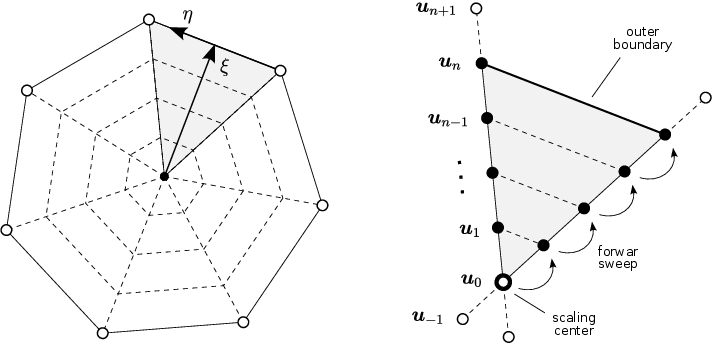}
        \caption{Schematics of radial discretization and condensation procedure.}
        \label{fig:procedure}
    \end{figure}

    From this point forward, we use the Thomas algorithm~\cite{patankar2018numerical} to eliminate the sub-diagonal terms $\bm\vartheta_{i}$ for $i\in\lbrace 1,2,\cdots,n \rbrace$. It proceeds by defining the pivot matrices
    \begin{equation}
        \bm\varrho_i = \bm\vartheta_i\bm\varphi^{-1}_{i-1},
    \end{equation}
    followed by performing the following forward sweeps
    \begin{equation}
        \begin{gathered}
            \bm\vartheta_i       := \bm\vartheta_i         - \bm\varrho_i\bm\varphi_{i-1} = \bm{0}, \\
            \bm\varphi_i         := \bm\varphi_i           - \bm\varrho_i\bm\psi_{i-1}, \\
            \bm p_i^{t+\Delta t} := \bm p_i^{t+\Delta t}   - \bm\varrho_i\bm p_{i-1}^{t+\Delta t}.
        \end{gathered}
    \end{equation}
    Finally, by defining
    \begin{equation}
        \bm\varrho_{n+1} = \bm\vartheta_{n+1}\bm\varphi^{-1}_{n-1},
    \end{equation}
    and applying
    \begin{equation}
        \begin{gathered}
            \bm\vartheta_{n+1}       := \bm\vartheta_{n+1}     - \bm\varrho_{n+1}\bm\varphi_{n-1} = \bm{0}, \\
            \bm\varphi_{n+1}         := \bm\varphi_{n+1}       - \bm\varrho_{n+1}\bm\psi_{n-1}, \\
            \bm p_{n+1}^{t+\Delta t} := -\bm\varrho_{n+1}\bm p_{n-1}^{t+\Delta t},
        \end{gathered}
    \end{equation}
    the last two equations are decoupled from the rest of the system, so that
    \begin{equation} \left[
        \begin{array}{ll}
            \bm\varphi_{n}  &\bm\psi_{n}  \\
            \bm\varphi_{n+1}&\bm\psi_{n+1}
        \end{array} \right] \left[
        \begin{array}{l}
            \bm{u}_{n}^{t+\Delta t} \\ \bm{u}_{n+1}^{t+\Delta t}
        \end{array} \right] = \left[
        \begin{array}{l}
            \bm{p}_n^{t+\Delta t} \\ \bm{p}_{n+1}^{t+\Delta t}
        \end{array} \right] + \left[
        \begin{array}{l}
            \bm{0} \\ \bm{f}^{t+\Delta t}
        \end{array} \right].
    \end{equation}
    The first equation gives
    \begin{equation}
        \bm{u}_{n+1}^{t+\Delta t} = \bm\psi_{n}^{-1}(\bm p_n^{t+\Delta t} - \bm\varphi_{n}\bm{u}_{n}^{t+\Delta t}).
    \end{equation}
    Substituting it into the second one yields
    \begin{equation}
        (\bm\varphi_{n+1} - \bm\psi_{n+1}\bm\psi_{n}^{-1}\bm\varphi_{n})\bm{u}_{n}^{t+\Delta t}  = \bm p_{n+1}^{t+\Delta t} - \bm\psi_{n+1}\bm\psi_{n}^{-1}\bm p_n^{t+\Delta t} + \bm f^{t+\Delta t},
    \end{equation}
    or, equivalently,
    \begin{equation}
        \bm{S} \bm{u}_n^{t+\Delta t} = \hat{\bm f}^{t+\Delta t},
    \end{equation}
    where
    \begin{equation}
        \bm{S} = \bm\varphi_{n+1} - \bm\psi_{n+1}\bm\psi_{n}^{-1}\bm\varphi_{n}
    \end{equation}
    is the dynamic stiffness matrix, and
    \begin{equation}
        \hat{\bm f}^{t+\Delta t} = \bm p_{n+1}^{t+\Delta t} - \bm\psi_{n+1}\bm\psi_{n}^{-1}\bm p_n^{t+\Delta t} + \bm{f}^{t+\Delta t}
    \end{equation}
    is the effective source vector. Recalling that $\bm{u}_n$ is in fact the nodal solution field, the global algebraic system for the entire domain can then be assembled and solved. As a result, the equation system in \eqref{eq:system}, which contains the internal information of the elements, does not need to be assembled explicitly. Instead, its effects are computed on the fly, keeping the dimensions of the global system unchanged. The schematic of this condensation procedure is depicted in Figure~\ref{fig:procedure}.
    
    Following the solution of the global system at time $t+\Delta t$ and finding the nodal field, the internal unknowns $\bm{u}_{i}$ can be computed using the backward sweep
    \begin{equation}
        \bm{u}_i^{t+\Delta t} = \bm\varphi_i^{-1}(\bm p_i^{t+\Delta t} - \bm\psi_i\bm{u}_{i+1}^{t+\Delta t}).
    \end{equation}
    Additionally, its time derivatives $\partial_t\bm{u}_i^{t+\Delta t}$ and $\partial^2_t\bm{u}_i^{t+\Delta t}$ are computed using~\eqref{eq:Newmark} and the solution proceeds to the next time step. Note that in case of equal time stepping, the stiffness matrix need be computed once, and only the forward and backward sweep for the residual source vectors and internal auxiliary unknowns must be performed.
    
    \section{Numerical implementation} \label{sec:numerical}
\begin{figure}
        \centering
        \includegraphics[scale=1.0]{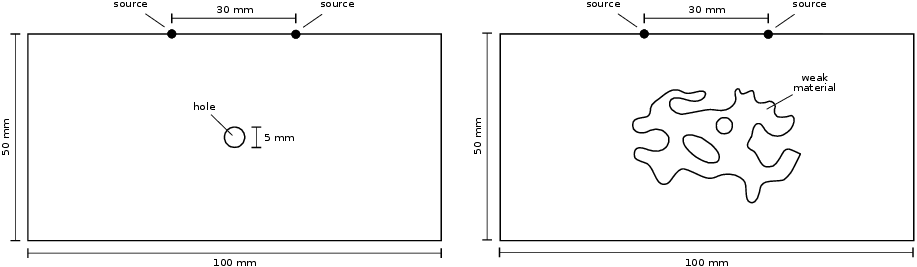}
        \caption{Configuration of the problems: rectangular specimen with a circular hole (left) and with a randomly shaped weak inclusion (right).}
        \label{fig:config}
    \end{figure}

    We test two rectangular specimens with the configurations shown in Figure~\ref{fig:config}. The first one contains a circular hole at the center and the second one a randomly-shaped weak inclusion. The density is $\rho=2700$ kg/m$^3$ everywhere, except for the weak part, which has half the density of the original material.
    
\begin{figure}
        \centering
        \includegraphics[scale=1.0]{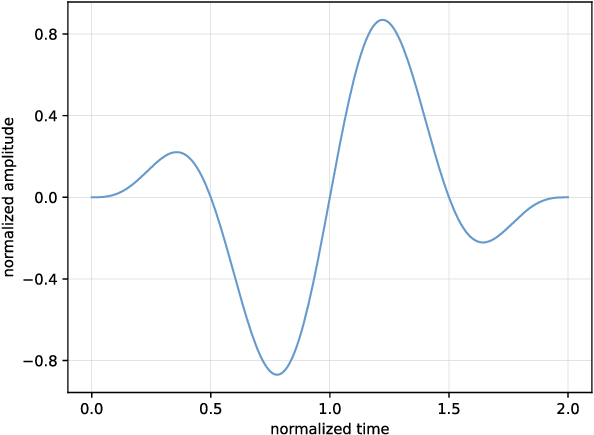}
        \caption{Sine burst impulse: curve of normalized amplitude $a(t)/a_0$ versus normalized time $tf_0$.}
        \label{fig:burst}
    \end{figure}

    Two sources with a pitch of 30 mm, centrally located on the top surface of the specimens are used to induce two forward wavefields. For this purpose, we use the sine burst
    \begin{equation}
        a(t) = 
        \begin{cases}
            a_0 \sin(\omega t) \sin^2({\omega t}/{2n_c}), & 0 \leq t \leq {n_c}/{f_0} \\
            0, & \text{otherwise}
        \end{cases}
    \end{equation}
    where $\omega = 2\pi f_0$ and $n_c$ is the number of cycles. Note that this signal is the commonly used impulse type in nondestructive testing. Similar to the work of Herrmann et al.~\cite{herrmann2023use}, we use a 2-cycle excitation with an amplitude of $a_0=10^{12}$ N/m$^2$ and a central frequency of $f_0=500$ kHz. The excitation signal curve is shown in Figure~\ref{fig:burst}. Note that the wave speed in the medium is taken to be $c=6000$ m/s, and we use a total time duration of 40 $\mu$s with 0.2 $\mu$s time steps.
    
    The measurements are made using equally spaced transducers with a pitch of 1 mm located along the four edges. The back-propagating wavefield is then induced by injecting the composite residual recorded at these points.
    
\begin{figure}
        \centering
        \includegraphics[scale=1.0]{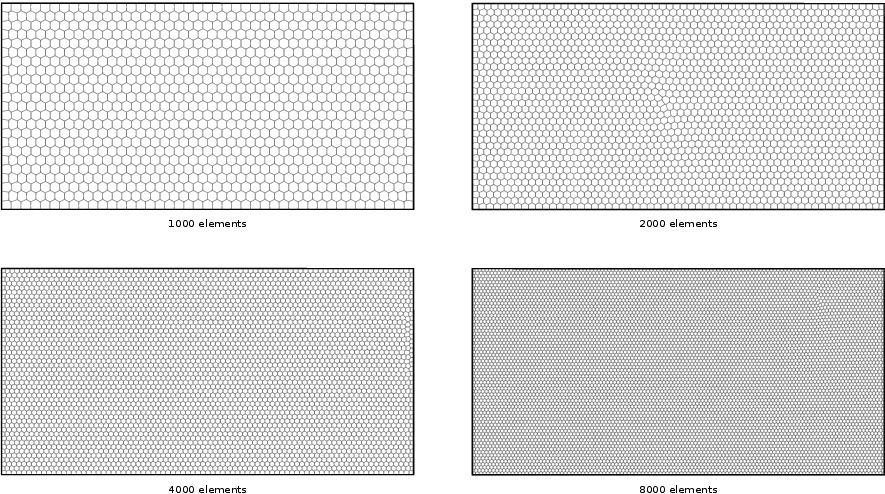}
        \caption{Centroidal Voronoi tessellations over the domain of problem with 1000, 2000, 4000, and 8000 elements.}
        \label{fig:mesh}
    \end{figure}

    Four meshes, consisting of 1000, 2000, 4000, and 8000 polygons, are defined (see Figure~\ref{fig:mesh}). To efficiently harness the strength of the scaled boundary finite element method, we use centroidal Voronoi tessellations, which provides the most possible symmetry for the polygons.
    
    The inversions are performed three times for each mesh. The first involves using the low-frequency expansion, while the second and third apply radial discretization with 5 and 7 finite difference intervals, respectively. Note that using fewer than 5 intervals can, in some cases, lead to getting stuck in local minima, which indicates a lack of accuracy in the simulated wavefields.
    
    It is worth noting that the formulation is implemented in an in-house C++ program. All computations are carried out using 12 threads on a laptop running Arch Linux 6.12.4, equipped with an AMD Ryzen 7 PRO 6850U CPU clocked at 2.7 GHz. We use the limited-memory BFGS algorithm~\cite{liu1989limited} provided by the NLopt nonlinear optimization library~\cite{NLopt2007} to move downhill on the cost function hypersurface. We begin all optimization problems by setting the constant density of $\rho_0=2700$ throughout the domain and applying a box constraint to the algorithm to ensure that the relative density $\rho/\rho_0$ remains above $10^{-6}$ and below 2. The lower bound prevents instability in the wave equation, while the upper bound prevents overshoot and ensures that the solution does not fall into unrealistic local minima. In addition, due to the placement of integration points in the scaled boundary finite element method, which are located on element boundaries, the gradient of the wavefields exhibits jumps between adjacent elements. Therefore, we omit the second term in the Fr\'echet kernel, as it involves the spatial gradients of the forward and adjoint wavefields. Consequently, we can define the density field at the nodal points and directly use the solution field of the forward and adjoint problems, and avoid an additional loop over elements. This approach significantly reduces the computational cost per epoch. In addition, this strategy enables efficient memory management by storing the forward wavefield generated by each source term and computing the Fr\'echet kernel on the fly while solving the adjoint problem.
    
    \subsection{Specimen with hole}
    First, we analyze the specimen containing a circular hole at its center. The synthetic data for the problem is generated using a boundary-conforming finite element model with 155,818 bilinear quadrilateral elements.
    
\begin{figure}
        \centering
        \includegraphics[scale=1.0]{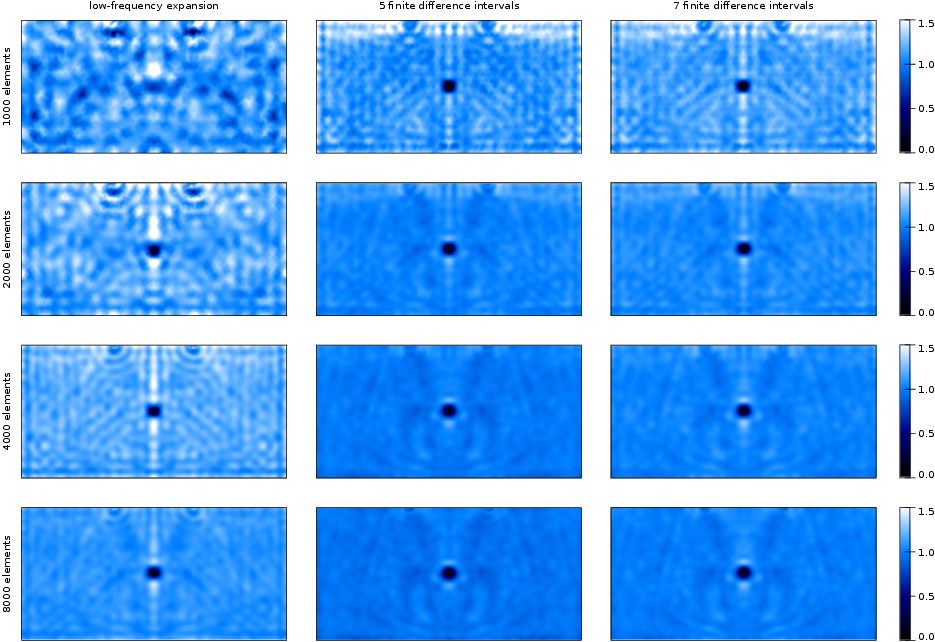}
        \caption{Normalized reconstructed density under different assumptions.}
        \label{fig:hole-comparison}
    \end{figure}

    The reconstructed density distribution of the specimen for different cases is presented in Figure~\ref{fig:hole-comparison}. The figure contains four rows, each designated for one of the meshes, and three columns for the type of approximation. The rows correspond to meshes with 1000, 2000, 4000, and 8000 elements, respectively. The first column represents the low-frequency expansion approximation, while the second and third columns represent radial discretization with 5 and 7 finite difference intervals, respectively. Elaborating on the first row, the low-frequency expansion is clearly unable to recover a meaningful representation of the real model, while radial discretization, although providing noisy distributions, captures the hole with either 5 or 7 finite difference intervals. Proceeding to the second row, the low-frequency expansion detects the hole, but the distribution appears more scattered and noisy compared to those of the radial discretization with the coarsest mesh. Both radial discretizations with 5 and 7 finite difference intervals provide fairly acceptable distributions. The trend of the third and fourth rows are similar to that of the second row, with the exception that as the number of elements increases, the low-frequency expansion solution becomes less noisy. However, the quality of the image reproduced by the finest mesh with 8000 elements is still lower than that reconstructed by the radial discretizations with 5 intervals using 2000 elements.
    
    Comparing the average wall-clock time per epoch for different cases, as presented in Table~\ref{tab:hole}, the latter case---using 2000 elements in the radial discretization with 5 intervals---took 2.39 seconds per epoch on average, while the low-frequency expansion solution with 8000 elements took 2.70 seconds. Regarding the reproduced images and computational costs, the new formulation clearly outperforms the existing method.
    
\begin{figure}
        \centering
        \includegraphics[scale=0.55]{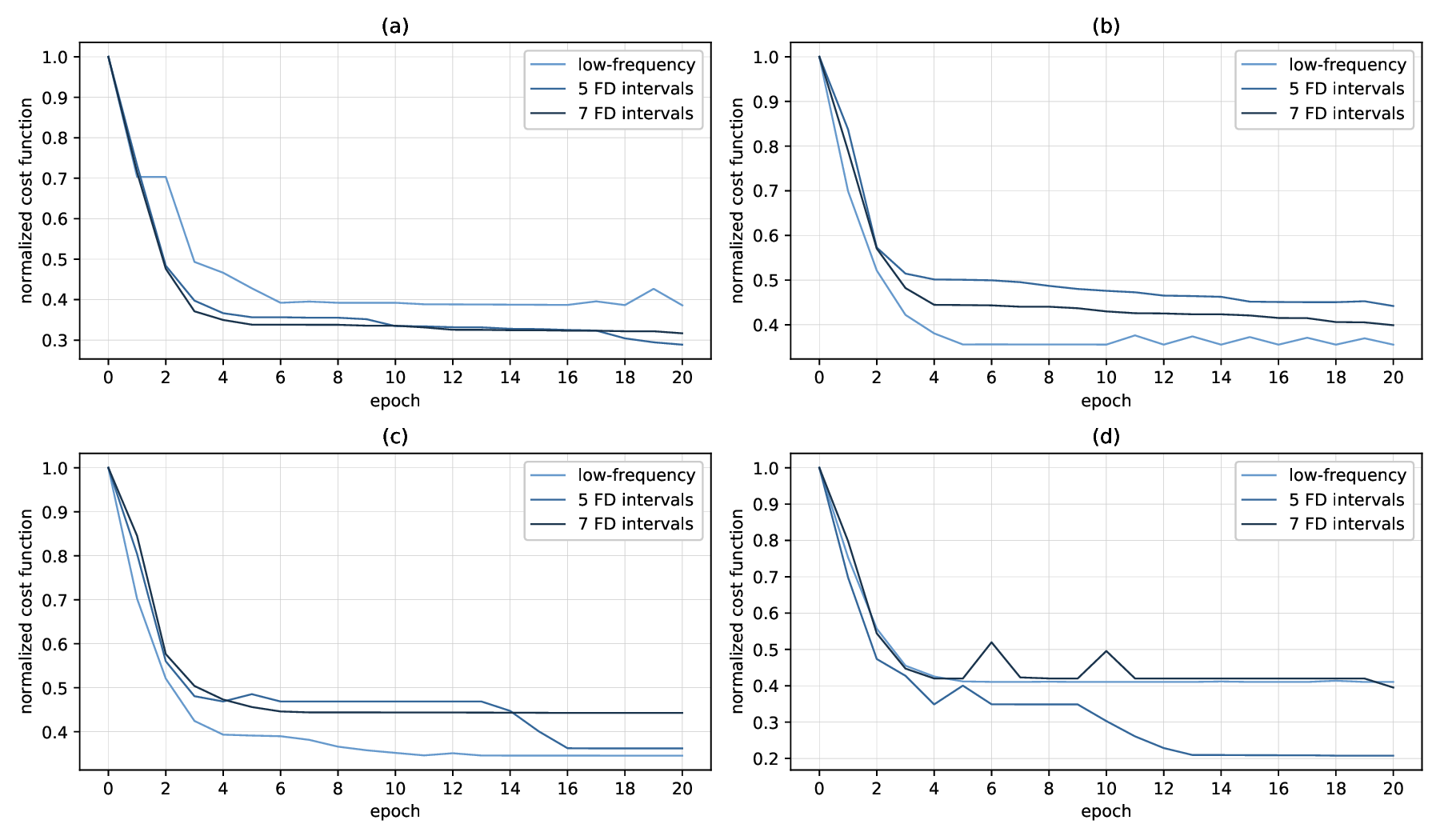}
        \caption{Normalized cost function over the course of epochs: (a) 1000 elements, (b) 2000 elements, (c) 4000 elements, and (d) 8000 elements.}
        \label{fig:hole-epochs}
    \end{figure}

\begin{figure}
        \centering
        \includegraphics[scale=0.75]{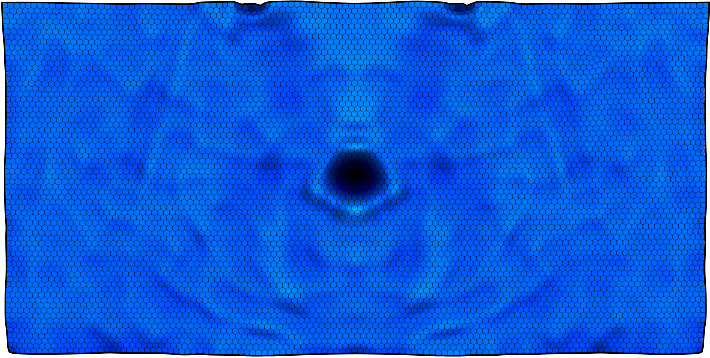}
        \caption{Three-dimensional representation of the normalized reconstructed density.}
        \label{fig:hole-3D}
    \end{figure}

    The normalized cost function for all cases over the course of epochs is presented in Figure~\ref{fig:hole-epochs}. Note that each curve is normalized by the value of its corresponding cost function at the initial epoch. A three-dimensional representation, with density depicted as the third dimension, for the most extreme case---finest mesh with 8000 elements using radial discretization with 7 finite difference intervals---is also presented in Figure~\ref{fig:hole-3D} to provide a better understanding of the reconstructed distribution.
    
    \begin{table}
        \centering
        \caption{Average wall-clock time per epoch for different cases in seconds.}
        \begin{tabular}{lccc}
            \hline        & low-frequency & 5 FD intervals & 7 FD intervals \\
            \hline
            1000 elements & 0.23          & 0.95           & 1.32           \\
            2000 elements & 0.46          & \textbf{2.39}  & 3.21           \\
            4000 elements & 1.27          & 5.49           & 6.82           \\
            8000 elements & \textbf{2.70} & 11.19          & 13.92          \\
            \hline
        \end{tabular}
        \label{tab:hole}
    \end{table}
    
    \subsection{specimen with weak inclusion}
    Next, we proceed to reconstruct the density distribution of the specimen containing a randomly-shaped weak inclusion. This relatively complicated shape is used to evaluate the capability of the presented formulation in a more challenging scenario. The weak region is assumed to have half the density of the original material to test whether the model can capture the two isolated regions within this weak part. Again, we use a boundary conforming finite element model to generate the synthetic data. It comprises 156,497 bilinear quadrilateral finite elements.
    
\begin{figure}
        \centering
        \includegraphics[scale=1.0]{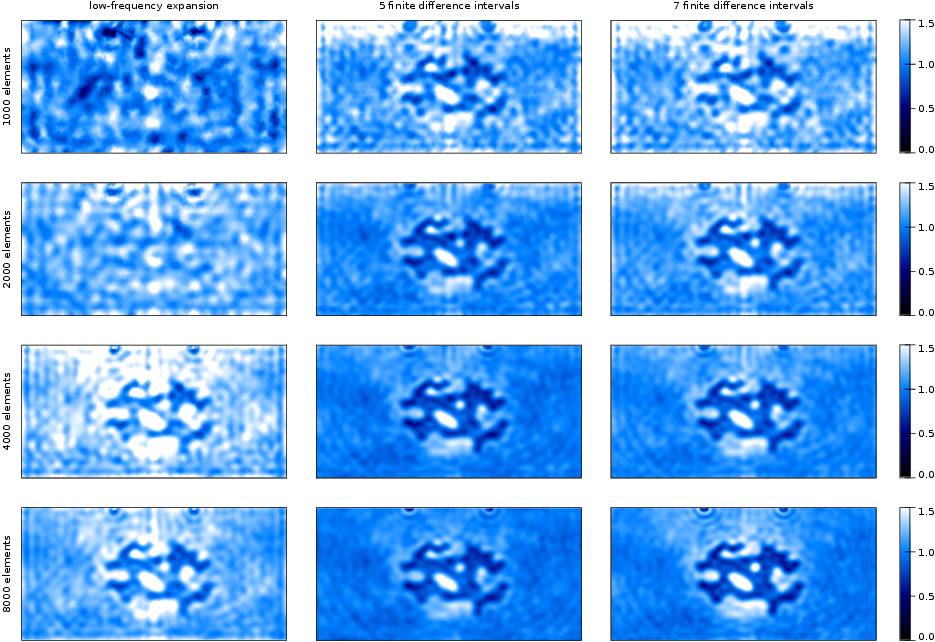}
        \caption{Normalized reconstructed density under different assumptions.}
        \label{fig:random-comparison}
    \end{figure}

    Figure~\ref{fig:random-comparison} presents the normalized reconstructed density of the specimen under different assumptions. The arrangement of the contours is similar to that of the previous specimen. A closer look at the first row shows that the low-frequency expansion fails to detect a meaningful shape. In contrast, the reconstructed densities obtained using radial discretization contain a low-quality image of the inclusion. The quality of the images produced by radial discretization improves significantly in the second row, while the low-frequency expansion continues to struggle in detecting the weak inclusion. From the third row onward, the inclusion becomes visible in the low-frequency expansion, but the images remain noisy. It requires more elements to eliminate these artifacts. On the other hand, radial discretization, whether with 5 or 7 finite difference intervals, produces sharp images.
    
    Table~\ref{tab:random} presents the average wall-clock time per epoch for different cases in seconds. Based on the data provided, with 2000 elements and 5 finite difference intervals, the new formulation achieved an acceptable image in just 2.52 seconds on average per epoch, whereas the simplified model, despite using 8000 elements, took longer (2.84 seconds on average per epoch) and resulted in a noisier image. This demonstrates that the new formulation outperforms the simplified model as it delivers higher-quality images in less time.
    
\begin{figure}
        \centering
        \includegraphics[scale=0.5]{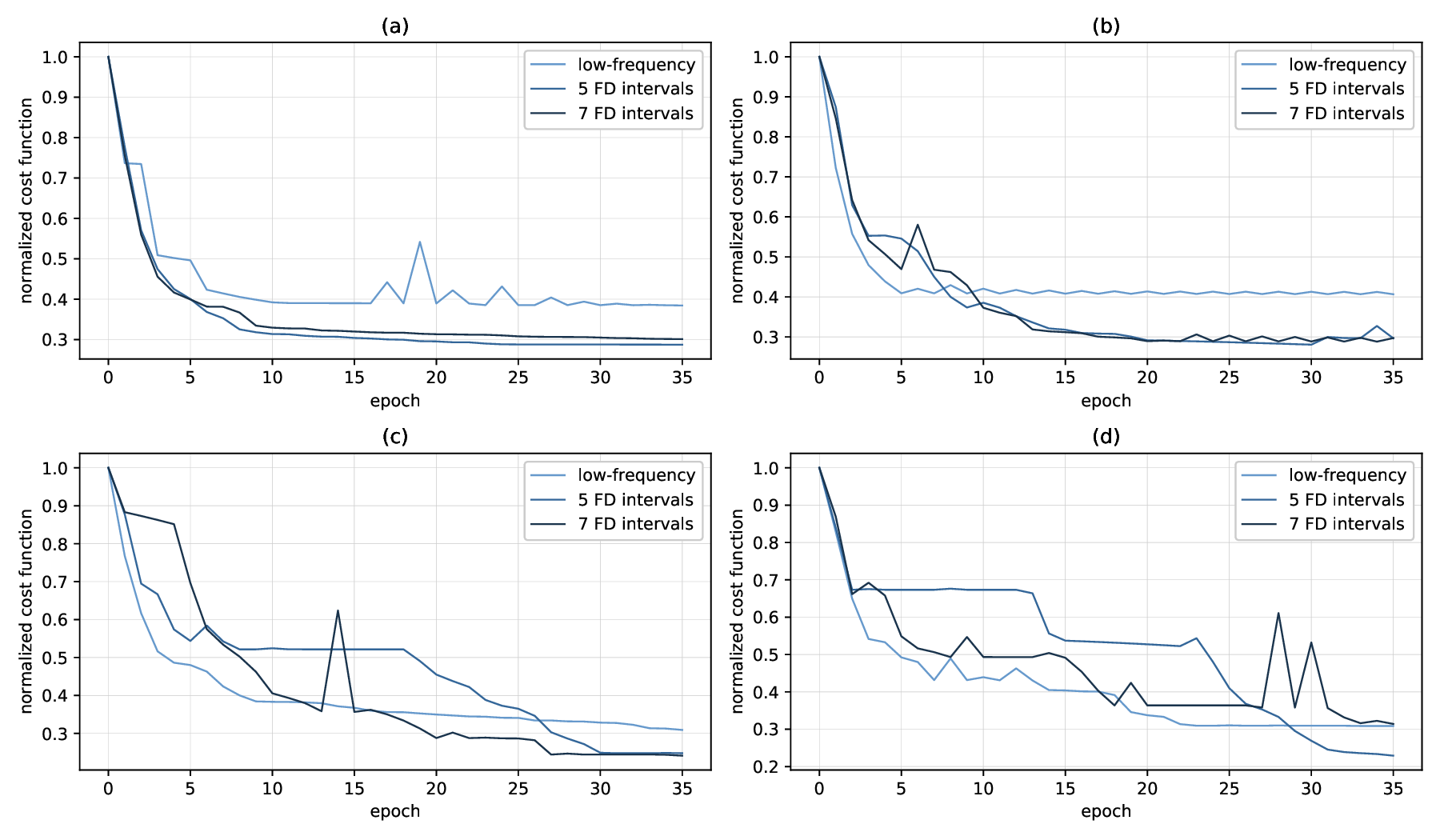}
        \caption{Normalized cost function over the course of epochs: (a) 1000 elements, (b) 2000 elements, (c) 4000 elements, and (d) 8000 elements.}
        \label{fig:random-epochs}
    \end{figure}

\begin{figure}
        \centering
        \includegraphics[scale=0.75]{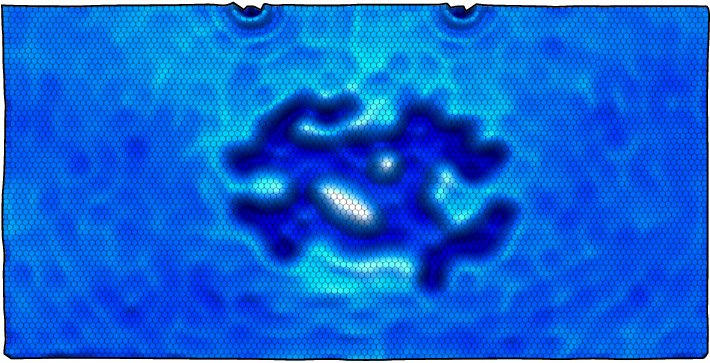}
        \caption{Three-dimensional representation of the normalized reconstructed density.}
        \label{fig:random-3D}
    \end{figure}

    Convergence trends are illustrated through the normalized cost function curves in Figure~\ref{fig:random-epochs}. It is evident that the minimization problem is more challenging than the previous example. Note that the absolute values of the cost function for the low-frequency expansion cases are far than those for radial discretization. Hence, we normalize each curve by its first value to visually represent the convergence trends. Finally, Figure~\ref{fig:random-3D} presents the three-dimensional representation of the case with 8000 elements and radial discretization using 7 finite difference intervals. The captured details are impressive.
    
    \begin{table}
        \centering
        \caption{Average wall-clock time per epoch for different cases in seconds.}
        \begin{tabular}{lccc}
            \hline        & low-frequency & 5 FD intervals & 7 FD intervals \\
            \hline
            1000 elements & 0.25          & 0.91           & 1.39           \\
            2000 elements & 0.46          & \textbf{2.52}  & 3.28           \\
            4000 elements & 1.51          & 5.41           & 7.16           \\
            8000 elements & \textbf{2.84} & 11.54          & 14.13          \\
            \hline
        \end{tabular}
        \label{tab:random}
    \end{table}
    
    \section{Conclusion} \label{sec:conclusion}
    We addressed the time-domain full-waveform inversion using the adjoint method and derived the scaled boundary semi-weak form of the scalar wave equation in heterogeneous media through the Galerkin method. We introduced a radial discretization scheme to solve the transient scalar wave equation with spatially varying density within the scaled boundary finite element framework. Two case studies---one with a specimen containing a hole and another with a weak inclusion---were used for demonstration. We generated synthetic data for the flawed samples using boundary-conforming finite element meshes to avoid committing inverse crime~\cite{wirgin2004inverse}. We utilized the limited-memory BFGS algorithm to efficiently determine the direction of maximum decrease in the cost function hypersurface. The numerical results demonstrated the effectiveness of our method in providing high-quality reconstructions even for complicated heterogeneous media. We also introduced an algorithmic condensation procedure to compute dynamic stiffness matrices on the fly, enabling a two-level hierarchical structure for the optimization problem. The use of local problems and a single coupled global system allows for efficient distribution of computational tasks. It significantly enhances the parallelizability of the solver and makes it highly scalable in high-performance computing environments. Accordingly, we showed that with a straightforward OpenMP implementation using 12 threads on a personal laptop, the new formulation outperforms the existing approach in both the quality of reconstructed images and computation time.
    
    \section{Authorship contribution statement}
    \textbf{Alireza Daneshyar:} Writing -- original draft, Conceptualization, Methodology, Visualization, Validation. 
    \textbf{Stefan Kollmannsberger:} Writing -- review \& editing, Conceptualization, Methodology, Supervision, Resources, Project administration.
    
    \section{Declaration of competing interest}
    The authors declare that they have no known competing financial interests or personal relationships that could have appeared to influence the work reported in this paper.
    
    \section{Data Availability}
    The implementation of the model is available at https://doi.org/10.5281/zenodo.14551328. This repository includes all code in C++ and the drivers to run all the examples.
    
    \section{Acknowledgement}
    The authors would like to acknowledge the financial support of the Alexander von Humboldt Foundation, Germany.
        
\bibliographystyle{ieeetr}
\bibliography{preprint}

\end{document}